%
%
%
%
\hsize=5in
\baselineskip=12pt
\vsize=20cm
\parindent=10pt
\pretolerance=40
\predisplaypenalty=0
\displaywidowpenalty=0
\finalhyphendemerits=0
\hfuzz=2pt
\frenchspacing
\footline={\ifnum\pageno=1\else\hfil\tenrm\number\pageno\hfil\fi}
%
%
\input amssym.def
\font\ninerm=cmr9
\font\ninebf=cmbx9
\font\ninei=cmmi9
\skewchar\ninei='177
\font\ninesy=cmsy9
\skewchar\ninesy='60
\font\nineit=cmti9
\def\reffonts{\baselineskip=0.9\baselineskip
	\textfont0=\ninerm
	\def\rm{\fam0\ninerm}%
	\textfont1=\ninei
	\textfont2=\ninesy
	\textfont\bffam=\ninebf
	\def\bf{\fam\bffam\ninebf}%
	\def\it{\nineit}%
\rm
	}
%
%
\def\frontmatter{\vbox{}\vskip1cm\bgroup
	\leftskip=0pt plus1fil\rightskip=0pt plus1fil
	\parindent=0pt
	\parfillskip=0pt
	\pretolerance=10000
	}
\def\endfrontmatter{\egroup\bigskip}
\def\title#1{{\font\titlef=cmbx12
  \titlef#1\par}}
\def\author#1{\bigskip#1\par}

\def\thanks#1{\footnote{}{\reffonts\rm\noindent#1\hfil}}
\def\section#1\par{\ifdim\lastskip<\bigskipamount\removelastskip\fi
	\penalty-250\bigskip
	\vbox{\leftskip=0pt plus1fil\rightskip=0pt plus1fil
	\parindent=0pt
	\parfillskip=0pt
  \pretolerance=10000{\bf#1}}\nobreak\medskip
	}
\def\proclaim#1. {\medbreak\bgroup{\noindent\bf#1.}\ \it}
\def\endproclaim{\egroup
	\ifdim\lastskip<\medskipamount\removelastskip\medskip\fi}
\newdimen\itemsize
\def\setitemsize#1 {{\setbox0\hbox{#1\ }
	\global\itemsize=\wd0}}
\def\item#1 #2\par{\ifdim\lastskip<\smallskipamount\removelastskip\smallskip\fi
	{\leftskip=\itemsize
	\noindent\hskip-\leftskip
	\hbox to\leftskip{\hfil\rm#1\ }#2\par}\smallskip}
\def\Proof#1. {\ifdim\lastskip<\medskipamount\removelastskip\medskip\fi
	{\noindent\it Proof\if\space#1\space\else\ \fi#1.}\ }
\def\endproof{\hfill\hbox{}\quad\hbox{}\hfill\llap{$\square$}\medskip}
\def\Remark. {\ifdim\lastskip<\medskipamount\removelastskip\medskip\fi
        {\noindent\bf Remark. }}
\def\endremark{\medskip}
%
%
\def\emph#1{{\it #1}\/}
\def\text#1{\hbox{#1}}
\def\mathrm#1{{\rm #1}}
%
%
\newcount\citation
\newtoks\citetoks
\def\citedef#1\endcitedef{\citetoks={#1\endcitedef}}
\def\endcitedef#1\endcitedef{}
\def\citenum#1{\citation=0\def\curcite{#1}%
	\expandafter\checkendcite\the\citetoks}
\def\checkendcite#1{\ifx\endcitedef#1?\else
	\expandafter\lookcite\expandafter#1\fi}
\def\lookcite#1 {\advance\citation by1\def\auxcite{#1}%
	\ifx\auxcite\curcite\the\citation\expandafter\endcitedef\else
	\expandafter\checkendcite\fi}
\def\cite#1{\makecite#1,\cite}
\def\makecite#1,#2{[\citenum{#1}\ifx\cite#2]\else\expandafter\clearcite\expandafter#2\fi}
\def\clearcite#1,\cite{, #1]}
%
%
\def\references{\section References\par
	\bgroup
	\parindent=0pt
	\reffonts
	\rm
	\frenchspacing
	\setbox0\hbox{99. }\leftskip=\wd0
	}
\def\endreferences{\egroup}
\newtoks\authtoks
\newif\iffirstauth
\def\checkendauth#1{\ifx\auth#1%
    \iffirstauth\the\authtoks
    \else{} and \the\authtoks\fi,%
  \else\iffirstauth\the\authtoks\firstauthfalse
    \else, \the\authtoks\fi
    \expandafter\nextauth\expandafter#1\fi
	}
\def\nextauth#1,#2;{\authtoks={#1 #2}\checkendauth}
\def\auth#1{\nextauth#1;\auth}
\newif\ifinbook
\newif\ifbookref
\def\nextref#1 {\par\hskip-\leftskip
	\hbox to\leftskip{\hfil\citenum{#1}.\ }%
	\initnextref}
\def\initnextref{\bookreffalse\inbookfalse\firstauthtrue\ignorespaces}
\def\paper#1{{\it#1},}
\def\InBook#1{\inbooktrue in ``#1",}
\def\book#1{\bookreftrue{\it#1},}
\def\journal#1{#1\ifinbook,\fi}
\def\BkSer#1{#1,}
\def\Vol#1{{\bf#1}}
\def\BkVol#1{Vol. #1,}
\def\publisher#1{#1,}
\def\Year#1{\ifbookref #1.\else\ifinbook #1,\else(#1)\fi\fi}
\def\Pages#1{\makepages#1.}
\long\def\makepages#1-#2.#3{\ifinbook pp. \fi#1--#2\ifx\par#3.\fi#3}
%
%
\newsymbol\square 1003
\newsymbol\rightrightarrows 1313
\let\Rar\Rightarrow
\let\Lrar\Leftrightarrow
\let\hrar\hookrightarrow
\let\ot\otimes
\let\sbs\subset
\def\co#1{^{\mkern1mu\mathrm{co}\mkern2mu#1}}
\def\cop{^{\mathrm{cop}}}
\def\id{\mathrm{id}}
\def\lco#1{\vphantom{H}\co{#1}\mkern-1mu}
\def\op{^{\mathrm{op}}}
\def\opcop{^{\mathrm{op,cop}}}
\def\triv{_{\mathrm{triv}}}
\def\defop#1#2{\def#1{\mathop{\rm #2}\nolimits}}
\defop\Coker{Coker}
\defop\Com{Com}
\defop\Cotor{Cotor}
\defop\Ext{Ext}
\defop\HOM{HOM}
\defop\Hom{Hom}
\defop\Img{Im}
\defop\Ker{Ker}
\defop\Tor{Tor}
\def\limdir{\mathop{\vtop{\offinterlineskip\halign{##\hskip0pt\cr\rm lim\cr
	\noalign{\vskip1pt}
	$\scriptstyle\mathord-\mskip-10mu plus1fil
	\mathord-\mskip-10mu plus1fil
	\mathord\rightarrow$\cr}}}}
\def\liminv{\mathop{\vtop{\offinterlineskip\halign{##\hskip0pt\cr\rm lim\cr
	\noalign{\vskip1pt}
	$\scriptstyle\mathord\leftarrow
	\mskip-10mu plus1fil\mathord-\mskip-10mu plus1fil\mathord-$\cr}}}}
\def\mapr#1{{}\mathrel{\smash{\mathop{\longrightarrow}\limits^{#1}}}{}}
\def\lmapr#1#2{{}\mathrel{\smash{\mathop{\count0=#1
  \loop
    \ifnum\count0>0
    \advance\count0 by-1\smash{\mathord-}\mkern-4mu
  \repeat
  \mathord\rightarrow}\limits^{#2}}}{}}
\def\mapd#1#2{\llap{$\vcenter{\hbox{$\scriptstyle{#1}$}}$}\big\downarrow
  \rlap{$\vcenter{\hbox{$\scriptstyle{#2}$}}$}}
\def\lmapd#1#2#3{\llap{$\vcenter{\hbox{$\scriptstyle{#2}$}}$}
  \left\downarrow\vcenter to#1pt{}\right.\mskip-4mu
  \rlap{$\vcenter{\hbox{$\scriptstyle{#3}$}}$}}
\def\lmapdr#1#2{\left\downarrow\vcenter to#1pt{}\right.\mskip-4mu
  \rlap{$\vcenter{\hbox{$\scriptstyle{#2}$}}$}}
\def\diagram#1{\vbox{\halign{&\hfil$##$\hfil\cr #1}}}
\let\al\alpha
\let\be\beta
\let\ep\varepsilon
\let\la\lambda
\let\ph\varphi
\let\De\Delta
\let\Si\Sigma
\def\calA{{\cal A}}
\def\calB{{\cal B}}
\def\calF{{\cal F}}
\def\calL{{\cal L}}
\def\calM{{\cal M}}
\def\calS{{\cal S}}
\def\calT{{\cal T}}
\def\0{_{(0)}}
\def\1{_{(1)}}
\def\2{_{(2)}}
\def\AM{\vphantom{\calM}_A\mskip1mu\calM}
\def\AMH{\AM^H}
\def\ApM{\vphantom{\calM}_{A'}\mskip1mu\calM}
\def\AT{\vphantom{\calT}_A\calT}
\def\ATH{\AT^H}
\def\BM{\vphantom{\calM}_B\mskip1mu\calM}
\def\BHM{\vphantom{\calM}_B^H\calM}
\def\CM{\vphantom{\calM}^C\!\calM}
\def\CopM{\vphantom{\calM}^{C\cop}\!\calM}
\def\CpM{\vphantom{\calM}^{C'}\!\calM}
\def\CT{\vphantom{\calM}^C\calT}
\def\DM{\vphantom{\calM}^D\!\calM}
\def\DMH{\vphantom{\calM}^D\!\calM_H}
\def\HCM{\vphantom{\calM}_H^{\,C\!}\calM}
\def\HCT{\vphantom{\calM}_H^{\,C}\calT}
\def\HCopM{\vphantom{\calM}_{H\cop}^{\,\,C\cop\!\!}\calM}
\def\HcoM{\vphantom{\calM}^{\,H\!}\calM}
\def\HopDM{\vphantom{\calM}_{H\op}^{\hphantom{H}\,\,D\!}\calM}
\def\HM{\vphantom{\calM}_H\calM}
\def\HMB{\vphantom{\calM}^H\!\calM_B}
\def\HMC{\vphantom{\calM}_H\calM^C}

\def\HMH{\vphantom{\calM}_H\calM^H}
\def\HHM{\vphantom{\calM}_H^{\,H\!}\calM}
\def\HTC{\vphantom{\calT}_H\calT^C}
\def\MA{\calM_A}
\def\MAH{\calM_A^H}
\def\MAK{\calM_A^K}
\def\MB{\calM_B}
\def\MC{\calM^C}
\def\MCp{\calM^{C'}}
\def\MD{\calM^D}
\def\MHD{\calM_H^D}
\def\MHH{\calM_H^H}
\def\TA{\calT_A}
\def\TAH{\calT_A^H}
\def\TC{\calT^C}
\def\sqC{{}\mathbin{\square}_C}
\def\sqCpH{{}\mathbin{\square}_{C'}\!H}
\def\sqCH{\sqC\!H}
\def\sqCM{\sqC M}
\def\sqD{{}\mathbin{\square}_D}
\def\sqDH{\sqD H}

\citedef
Bia-HM63
Br99
Br-W
Cae-G04
Cl-PS77
Das-NR
Doi83
Doi92
Ei-M65
Fai
Gab62
Gom-N95
Gross
Ja
Kop93
Lin77
Ma91
Ma94
Ma-D92
Ma-W94
Nas-T96
Nas-TZ96
Ob77
Par-W91
Pro
Scha00
Scha-Schn05
Schn92
Schn93
Scho
Sk07
Sk10
Sk21
St
Sw
Tak77
Tak79
Tak94
Ulb90
\endcitedef

\frontmatter
\title{On Takeuchi's correspondence}
\author{Serge Skryabin}
\endfrontmatter

\bigskip
{\reffonts
\leftskip=20pt\rightskip=20pt
{\bf Abstract.}
In this paper we review the Takeuchi correspondence between right coideal 
subalgebras and left $H$-module factor coalgebras of a Hopf algebra with 
bijective antipode. We are especially interested to describe the situation 
when the faithfulness assumption in the conditions of flatness and coflatness 
is dropped. This study was motivated by several questions which remain 
open.\par
}

\bigskip
\section
Introduction

Let $H$ be a Hopf algebra over a field $k$. In his 1979 paper \cite{Tak79} 
Takeuchi considered two classes of objects: the \emph{right coideal subalgebras} 
$A\sbs H$, i.e., subalgebras such that $\De(A)\sbs A\ot H$ with respect to 
the comultiplication $\De$ in $H$, and the \emph{left $H$-module factor 
coalgebras} $C=H/I$ which are quotients of $H$ with $I$ being simultaneously a 
coideal and a left ideal of $H$. To a right coideal subalgebra $A$ there 
corresponds the left $H$-module factor coalgebra $H/HA^+$ where $A^+$ is the 
augmentation ideal of $A$. In the opposite direction, given a left $H$-module 
factor coalgebra $C$, all elements of $H$ invariant under the left coaction of 
$C$ on $H$ form a right coideal subalgebra $\lco CH$. Takeuchi proved several 
important facts under the assumptions of either faithful flatness or faithful 
coflatness. These results captured essential features of the correspondence 
between subgroups of an affine algebraic group and affine homogeneous spaces.

As was elucidated by Masuoka and Wigner \cite{Ma-W94, Th. 2.1}, there are even 
somewhat better conclusions in the case when \emph{the antipode of $H$ is bijective}. 
Bijectivity of the antipode $S:H\to H$ will be essential for all central 
results of our present paper. So we take this condition as a standing 
assumption for the whole paper which will not be repeated in statements.

The starting point for our discussions is the fundamental fact, due to 
Takeuchi and Masuoka, which will be reproduced in our present paper as

\proclaim
Theorem 1.7.
Takeuchi's correspondence gives a bijection between the set of right coideal 
subalgebras of $H$ over which $H$ is left (respectively, right) faithfully flat 
and the set of left $H$-module factor coalgebras of $H$ over which $H$ is left 
(respectively, right) faithfully coflat.
\endproclaim

Unfortunately, this fact in its full generality is rarely mentioned in 
literature. In \cite{Tak79} such a bijection was stated only for a Hopf 
algebra which is either commutative or cocommutative. The reason is that for 
a Hopf algebra whose antipode is not bijective one does not get the conclusion 
of Theorem 1.7, and so in this case one is forced to consider narrower classes 
of respective objects. Notably, Schneider \cite{Schn93} described a bijection 
between normal Hopf subalgebras and conormal Hopf factor algebras. Several 
other bijections of this kind were presented in \cite{Ma-D92} and \cite{Tak94}.

The ultimate conclusion for a Hopf algebra with bijective antipode was finally 
spelled out by Masuoka \cite{Ma94, Th. 1.11}. However, it was merely noted 
there that one has this conclusion by the Masuoka-Wigner Theorem. Strictly 
speaking, one needs also dual arguments. We will give more details on Theorem 
1.7, viewing it as an immediate extension of original Takeuchi's theorems 
aided by a few easy observations concerning switchings between the left and 
right sides.

Faithful flatness and coflatness are very strong conditions. Considered 
separately, the property of $H$ to be a projective $A$-module or a generator 
in the category of $A$-modules corresponds to the property of $H$ to be a 
cogenerator or injective in the category of comodules over the left $H$-module 
factor coalgebra $C$ corresponding to a right coideal subalgebra $A$. We 
review such conclusions in section 2 of this paper. They are simpler variants 
of results obtained by Doi \cite{Doi92} in a more general situation.

It seems to be less known that flatness of $H$ over its right coideal 
subalgebra $A$ can also be recognized by means of the corresponding factor 
coalgebra. According to the definition given in \cite{Scho, p.~122} and 
\cite{St, p.~225} the \emph{dominion} of $A$ in $H$ is the subalgebra of $H$ 
consisting of all its elements $h$ such that
$$
h\ot1=1\ot h\quad\text{in $H\ot_AH$}.
$$
We say that $A$ is a \emph{dominion subalgebra} if $A$ coincides with its own 
dominion in $H$. The relevance of this notion is explained by the fact that 
the right coideal dominion subalgebras of $H$ are precisely the right coideal 
subalgebras of the form $\lco CH$ for left $H$-module factor coalgebras $C$ of 
$H$. The bijection of Theorem 1.7 is now extended to larger sets as follows:

\proclaim
Theorem 4.6.
Takeuchi's correspondence gives a bijection between the set of right coideal 
dominion subalgebras of $H$ over which $H$ is left (respectively, right) flat 
and the set of left $H$-module factor coalgebras $C$ of $H$ with the property 
that each finite-dimensional left (respectively, right) $C$-comodule embeds 
as a $C$-subcomodule in a left (respectively, right) $H$-comodule.
\endproclaim

In the case when $H=k[G]$ is the algebra of regular functions on an affine 
algebraic group $G$ the extendibility property occurring in Theorem 4.6 
distinguishes \emph{observable subgroups}\/: a closed subgroup $K$ of $G$ is 
observable if each finite-dimensional rational $K$-module embeds as a 
$K$-submodule in a rational $G$-module. By an old result of Bia{\l}ynicki-Birula, 
Hochschild and Mostow \cite{Bia-HM63} this property characterizes subgroups 
such that the respective quotient $G/K$ is a quasi-affine algebraic variety 
(see Grosshans \cite{Gross} for more information). Moreover, it was shown in 
\cite{Bia-HM63} that it suffices to assume the extendibility of 1-dimensional 
representations. We are not able to obtain such a strengthening of Theorem 
4.6 as there are no general constructions of 1-dimensional comodules similar 
to the top exterior powers in the commutative theory.

As was shown by Schauenburg \cite{Scha00} a noncommutative Hopf algebra is not 
necessarily a flat module even over its Hopf subalgebras. An interesting 
feature we come across is that the comodule extendibility property of a factor 
coalgebra in Theorem 4.6 does ensure flatness of $H$ over the corresponding 
coideal subalgebra.

Even in the class of commutative Hopf algebras, it is quite possible for $H$ 
to be a flat module over its right coideal subalgebra $A$ but not faithfully 
flat. If $H=k[G]$ for an affine algebraic group $G$, this happens precisely 
when $A$ corresponds to a quasi-affine homogeneous space $G/K$ which is not 
affine. The dual situation is very much different, and we may ask

\proclaim
Question.
Is there any example of a left $H$-module factor coalgebra $C$ over 
which $H$ is left or right coflat but not faithfully coflat?
\endproclaim

By a result of Doi \cite{Doi83, p.~247} coflatness does imply faithful 
coflatness when $C$ is a factor Hopf algebra of $H$. This result gives a Hopf 
algebraic interpretation of the well-known characterization of affine 
quotients $G/K$ in terms of exactness property of the induction functor 
\cite{Cl-PS77}, \cite{Ob77}. If $H$ has cocommutative coradical, then $H$ is 
faithfully coflat over all $H$-module factor coalgebras by a theorem of 
Masuoka \cite{Ma91}. In section 7 of our present paper it will be shown that 
coflatness implies faithful coflatness under several different conditions on a 
left $H$-module factor coalgebra $C$ of $H$. Particularly useful is

\proclaim
Theorem 7.4.
Suppose that $H$ is left coflat over $C$ and right flat over $A=\lco CH$. Then 
$H$ is left faithfully coflat over $C,$ left faithfully flat over $A$.

If, in addition, either $C$ has a proper right coideal of finite codimension 
or all simple left $A$-modules are finite dimensional, then $H$ is also right
faithfully coflat over $C,$ right faithfully flat over $A$.

\endproclaim

It can happen that $H$ is flat over all right coideal subalgebras. Then 
Theorem 7.4 shows that for any left $H$-module factor coalgebra $C$ coflatness 
of $H$ over $C$ implies faithful coflatness. This argument works in the 
classical case of commutative Hopf algebras, and moreover in the PI case --- 
see Corollary 7.7.

The hypothesis of Theorem 7.4 demands flatness and coflatness on opposite 
sides. As compared with the flatness and coflatness on the same side, our 
assumptions in Theorem 7.4 lead to stronger conclusions and avoid the use of 
Galois condition (cf. \cite{Scha-Schn05, Prop. 4.5}).

Takeuchi's work in \cite{Tak79} made heavy use of the categories of relative 
Hopf modules. For example, $\MAH$ is the category of right-right $(H,A)$-Hopf 
modules. Its objects are right $A$-modules equipped with a right $H$-comodule 
structure such that the map $M\ot_kA\to M$ given by the action of $A$ on $M$ 
is $H$-colinear with respect to the tensor product of comodule structures on 
$M$ and $A$. There are several similar categories considered in \cite{Doi92}, 
\cite{Tak79} and other papers.

We denote by $\MA$ and $\AM$ the categories of right and left 
$A$-modules, by $\MC$ and $\CM$ the categories of right and left $C$-comodules. 
Consider the following conditions on a right coideal subalgebra $A$ and a left 
$H$-module factor coalgebra $C$:
$$
\openup1\jot
\eqalignno{
&\text{all objects of $\MAH$ are projective in $\MA$}&(0.1)\cr
&\text{all objects of $\HCM$ are injective in $\CM$}&(0.2)\cr
&\text{all objects of $\AMH$ are projective in $\AM$}&(0.3)\cr
&\text{all objects of $\HMC$ are injective in $\MC$}&(0.4)\cr
}
$$
If $H$ is right faithfully flat over $A$, then (0.1) holds in the case when 
$A$ is contained in the center of $H$ \cite{Tak79, Th. 5}, and also in the 
case when $A$ is a Hopf subalgebra \cite{Doi83, Th. 4}. If $H$ is right 
$A$-cocleft, then each object of $\MAH$ is even a free $A$-module 
\cite{Ma-D92, p.~3715}. Motivated by these results Masuoka asked whether each 
nonzero object of $\AMH$ is a projective generator in $\AM$ (this is 
equivalent to condition (0.3)) whenever $H$ is left faithfully flat over $A\,$ 
\cite{Ma94, Question 1.9}. 

We will see in Lemma 8.1 that $H$ is faithfully flat over $A$ on both sides 
whenever (0.1) holds. Here condition (0.1) can be replaced by (0.3). Therefore 
Masuoka's question implies another one: \emph{is faithful flatness of $H$ over 
$A$ a left-right symmetric property}\/? We are not able to answer the latter 
question either, although this symmetry has been known in the case when $A$ is 
a subbialgebra \cite{Schn92, Cor. 1.8}, and another special case can be 
deduced from the already mentioned Theorem 7.4 --- see Corollary 7.6. However, 
one could try to modify Masuoka's question by imposing two-sided faithful 
flatness. We will prove

\proclaim
Theorem 8.4.
Takeuchi's correspondence gives a bijection between the set of right coideal 
subalgebras of $H$ satisfying condition $(0.1)$ and the set of left $H$-module 
factor coalgebras of $H$ satisfying condition $(0.2)$.

There is a similar bijection with respect to conditions $(0.3)$ and $(0.4)$.
\endproclaim

\proclaim
Theorem 8.6.
Suppose that $H$ is left and right faithfully flat over its right coideal 
subalgebra $A,$ and let $C=H/HA^+$. If either $A$ is an algebra of finite weak 
global dimension or $C$ is a coalgebra of finite global dimension, then
properties $(0.1)$--$(0.4)$ are all satisfied.
\endproclaim

The \emph{global dimension} of a coalgebra $C$ is the homological dimension of 
the abelian categories $\MC$ and $\CM$. This notion was introduced by 
N\u{a}st\u{a}sescu, Torrecillas and Zhang \cite{Nas-TZ96}. By \cite{Nas-TZ96, 
Cor. 3} $\MC$ and $\CM$ have equal homological dimensions defined as the 
supremum of injective dimensions of comodules. The \emph{weak global 
dimension} of algebras is defined in terms of lengths of flat resolutions of 
modules. This dimension is known to be left-right symmetric.

It is not clear whether finiteness of global dimensions in Theorem 8.6 is 
necessary. Without this assumption we will prove in section 9 several results 
on vanishing of $\Ext$ and $\Tor$ which serve as approximation to the desired 
conclusion.

In the case when $H$ is left faithfully flat over its right coideal subalgebra 
$A$ one of Takeuchi's results in \cite{Tak79} gives an equivalence of 
categories $\MAH\approx\MC$. If $H$ is flat over $A$, but not faithfully flat, 
there is a similar equivalence with $\MAH$ replaced by a certain quotient 
category. This was shown in \cite{Sk10} under the assumption that the algebras 
$A$ and $H$ have artinian classical quotient rings. In section 10 of our 
present paper we will revisit this result and explain that the only thing 
needed here is the flatness assumption. We will also generalize the second 
Takeuchi equivalence $\HCM\approx\AM$. The results of section 10 provide 
examples of perfect localizations of algebras and perfect colocalizations of 
coalgebras.

I would like to thank Akira Masuoka for early comments and clarification of 
his approach to \cite{Ma94, Th. 1.11}.

\section
1. Initial results

We will work over an arbitrary ground field $k$. The subscript $k$ in the 
notation for the tensor product functor $\ot_k$ will be omitted. Let $H$ be a 
Hopf algebra over $k$ with the comultiplication $\De:H\to H\ot H$, the counit 
$\ep:H\to k$ and the antipode $S:H\to H$. Unless explicitly stated otherwise 
it will be assumed that \emph{$S$ is bijective}. Put $H^+=\Ker\ep$.

If $A$ is a right or left coideal of $H$, then 
$A^+=A\cap H^+$ is a coideal of $H$, whence $HA^+$ is simultaneously a coideal 
and a left ideal of $H$. The quotient $H/HA^+$ is then a left $H$-module 
factor coalgebra of $H$.

Suppose that $C$ is an arbitrary factor coalgebra of $H$. Let $\pi:H\to C$ 
be the canonical surjective homomorphism of coalgebras. Viewing $H$-comodules 
as $C$-comodules via $\pi$ we obtain canonical functors $\calM^H\to\MC$ and 
$\HcoM\to\CM$. They admit right adjoint functors $?\sqCH$ and $H\sqC{}?$ 
given by the cotensor products with $H$ over $C$ (see \cite{Tak77, Appendix 2} 
and \cite{Das-NR, Ch.  2} for details on cotensor products).

The image $1_C=\pi(1)$ of the identity element $1\in H$ is a \emph{distinguished 
grouplike} element of $C$. For each left or right $C$-comodule $V$ the subspace 
of coaction invariants is defined, respectively, by the rule
$$
\openup1\jot
\eqalign{
\lco CV=\{x\in V\mid\rho(x)&=1_C\ot x\}\cong k1_C\sqC V,\cr
V\co C=\{x\in V\mid\rho(x)&=x\ot 1_C\}\cong V\sqC k1_C
}
$$
where $\rho$ is the comodule structure map $V\to C\ot V$ or $V\to V\ot C$, 
respectively. In particular,
$$
\lco CH=\{h\in H\mid(\pi\ot\id)\De(h)=1_C\ot h\}
$$
is a right coideal, while $H\co C$ is a left coideal of $H$. 
The following fact is verified easily:

\proclaim
Lemma 1.1.
Let $A$  be a right coideal, $B$ a left coideal, $C$ a factor coalgebra of the 
Hopf algebra $H$. 
Then $A\sbs\lco CH$ $(B\sbs H\co C)$ if and only if $\pi(x)=\ep(x)1_C$ for all 
$x\in A\ (\text{respectively},\ x\in B$). Assuming that $1\in A$ and $1\in B,$ 
we thus have
$$
A\sbs\lco CH\Lrar A^+\sbs\Ker\pi,\qquad\qquad
B\sbs H\co C\Lrar B^+\sbs\Ker\pi.
$$
\endproclaim

Further on we will assume that $A$ is a \emph{right coideal subalgebra} of $H$ 
and $C$ a \emph{left $H$-module factor coalgebra} of $H$. The canonical surjection 
$\pi:H\to C$ is then left $H$-linear. This implies that the one-sided coideals 
$\co CH$ and $H\co C$ are subalgebras of $H$.

By Lemma 1.1 $A\sbs\lco CH$ if and only if $HA^+\sbs\Ker\pi$. If these 
inclusions hold, then there is a pair of adjoint functors
$$
\openup1\jot
\vcenter{\halign{\hfil$#$&\qquad\qquad\hfil$#$\cr
\Phi:\MAH\to\MC&\Psi:\MC\to\MAH\cr
M\mapsto M/MA^+,&V\mapsto V\sqCH\cr
}}\eqno(1.1)
$$
and another pair of adjoint functors
$$
\openup1\jot
\vcenter{\halign{\hfil$#$&\qquad\qquad\hfil$#$\cr
\Phi:\AM\to\HCM&\Psi:\HCM\to\AM\cr
V\mapsto H\ot_AV,&M\mapsto\lco CM.\cr
}}\eqno(1.2)
$$
Now we recall two main results of \cite{Tak79}. The statements we give here 
are slightly weaker than those in \cite{Tak79}:

\proclaim
Takeuchi's Theorem 1.
Suppose that $\,C=H/HA^+$ and $H$ is left faithfully flat over $A$. Then 
$\,A=\lco CH$ and $H$ is left faithfully coflat over $C$. Moreover, the two 
functors in $(1.1)$ are quasi-inverse equivalences.
\endproclaim

\proclaim
Takeuchi's Theorem 2.
Suppose that $\,A=\lco CH$ and $H$ is right faithfully coflat over $C$. Then 
$\,C=H/HA^+$ and $H$ is right faithfully flat over $A$. Moreover, the two 
functors in $(1.2)$ are quasi-inverse equivalences.
\endproclaim

In Takeuchi's Theorems bijectivity of the antipode $S$ is not required. In the 
case when $S$ is bijective we can apply these theorems to the Hopf algebras 
$H\op$ and $H\cop$ obtained from $H$ by changing either the multiplication or 
the comultiplication to the opposite one.

\proclaim
Corollary 1.2.
Suppose that $H$ is right faithfully flat over $A,$ and let $D=H/A^+H$. Then 
$A=\lco DH$ and $H$ is left faithfully coflat over $D$. The functors
$$
\openup1\jot
\vcenter{\halign{\hfil$#$&\qquad\qquad\hfil$#$\cr
\Phi:\AMH\to\MD&\Psi:\MD\to\AMH\cr
M\mapsto M/A^+M,&V\mapsto V\sqDH\cr
}}\eqno(1.3)
$$
are quasi-inverse equivalences.
\endproclaim

\Proof.
Viewing $A\op$ as a right coideal subalgebra of $H\op$, we may identify 
$\AMH$ with $\calM_{A\op}^{H\op}$, and the latter category is equivalent 
to $\MD$ by means of functors given in Takeuchi's Theorem 1.
\endproof

\proclaim
Corollary 1.3.
Suppose that $H$ is left faithfully coflat over $C,$ and let $B=H\co C$. 
Then $C=H/HB^+$ and $H$ is right faithfully flat over $B$. The functors
$$
\openup1\jot
\vcenter{\halign{\hfil$#$&\qquad\qquad\hfil$#$\cr
\Phi:\BM\to\HMC&\Psi:\HMC\to\BM\cr
V\mapsto H\ot_BV,&M\mapsto M\co C\cr
}}\eqno(1.4)
$$
are quasi-inverse equivalences.
\endproclaim

\Proof.
Here $C\cop$ is a left $H\cop$-module factor coalgebra of $H\cop$, and we get 
category equivalences $\,\HMC\approx\HCopM\approx\BM\,$ by Takeuchi's Theorem 2.
\endproof

Note that in the above corollaries $D$ is a right $H$-module factor coalgebra, 
while $B$ is a left coideal subalgebra of $H$. We have to compare $D$ and $B$ 
with $C$ and $A$ in Takeuchi's Theorems.

\proclaim
Lemma 1.4.
If $A$ is a right coideal subalgebra of $H,$ then $S(A)$ is a left coideal 
subalgebra and $A^+H=S(A)^+H$. If $B$ is a left coideal subalgebra of $H,$ 
then $S(B)$ is a right coideal subalgebra and $HB^+=HS(B)^+$.
\endproclaim

\Proof.
Since the antipode $S$ is an antiendomorphism of $H$ as an algebra and as a 
coalgebra, it sends subalgebras to subalgebras and left coideals to 
right coideals. Now $S(A)^+=S(A^+)$, and the equality $A^+H=S(A^+)H$ is 
observed in Koppinen's \cite{Kop93, Lemma 3.1}. The second assertion of Lemma 
1.4 is equivalent to the first one with $H$ replaced by $H\opcop$.
\endproof

\proclaim
Lemma 1.5.
If $I$ is a left ideal coideal of $H,$ then its full preimage $S^{-1}(I)$ 
under the antipode is a right ideal coideal and 
$\ H\co{H/I}=H\co{H/S^{-1}(I)}$. 

Similarly, if $J$ is a right ideal coideal of $H,$ then $S^{-1}(J)$ is a left 
ideal coideal and\quad$\lco{H/J}H=\lco{H/S^{-1}(J)}H$.
\endproclaim

\Proof.
If $x\in S^{-1}(I)$, then $S(x)\in I$. We get $xh\in S^{-1}(I)$ for all $h\in H$ 
since $S(xh)=S(h)S(x)\in I$, and $\De(x)\in S^{-1}(I)\ot H+H\ot S^{-1}(I)$ since
$$
(S\ot S)\De(x)=\De\op(Sx)\in I\ot H+H\ot I.
$$
Thus $S^{-1}(I)$ is a right ideal and a coideal of $H$. The set $H\co{H/I}$ 
consists of all elements $y\in H$ such that
$$
\De(y)-y\ot1\in H\ot I.\eqno(1.5)
$$
The set $H\co{H/S^{-1}(I)}$ is characterized by means of similar inclusions 
which can be also rewritten as
$$
(\id\ot S)\De(y)-y\ot1\in H\ot I.\eqno(1.6)
$$
Define linear transformations $\ph$ and $\psi$ of $H\ot H$ by the formulas
$$
\eqalign{
\ph(u\ot v)&=\De(u)\cdot(1\ot v)=\sum u\1\ot u\2v,\cr
\psi(u\ot v)&=(\id\ot S)\De(u)\cdot(1\ot v)=\sum u\1\ot S(u\2)v
}
$$
where $u,v\in H$. Note that $\ph=\psi^{-1}$ and $H\ot I$ is stable under both 
$\ph$ and $\psi$ since $I$ is a left ideal of $H$. Hence $\psi(H\ot I)=H\ot I$. 
Since
$$ 
\psi\bigl(\De(y)\bigr)=\psi\ph(y\ot1)=y\ot1,\qquad\psi(y\ot1)=(\id\ot S)\De(y),
$$
we see that for each $y\in H$ containment (1.5) is equivalent to (1.6).
\endproof

Bijectivity of the antipode is not required in Lemmas 1.4 and 1.5. If  $S$ is 
bijective, the equalities of Lemma 1.4 can be rewritten as
$$
A^+H=S(HA^+),\qquad HB^+=S(B^+H).\eqno(1.7)
$$
Lemma 1.5 can be reformulated as follows:

\proclaim
Lemma 1.6.
Let $C$ be a left $H$-module factor coalgebra, $D$ a right $H$-module factor 
coalgebra of $H$. Then
$$
\lco CH=S(\,H\co C),\qquad H\co D=S(\lco DH).
$$ 
\endproclaim

\Proof.
We have $C=H/I$ where $I$ is a left ideal coideal of $H$. By Lemma 1.5 
$H\co C$ consists of all elements $y\in H$ satisfying (1.6). Since $S$ is 
bijective, (1.6) is equivalent to the containment 
$\,(S\ot S)\De(y)-Sy\ot1\in H\ot I$. This can be rewritten as
$$
\De(Sy)-1\ot Sy\in I\ot H,
$$
which means precisely that $\,Sy\in\lco CH$.
\endproof

In our formulation of Masuoka's theorem \cite{Ma94, Th. 1.11} given below we 
treat both left side and right side conditions. We include a different proof 
which is linked very closely to Takeuchi's Theorems 1 and 2.

\proclaim
Theorem 1.7.
Takeuchi's correspondence gives a bijection between the set of right coideal 
subalgebras of $H$ over which $H$ is left (respectively, right) faithfully flat 
and the set of left $H$-module factor coalgebras of $H$ over which $H$ is left 
(respectively, right) faithfully coflat.
\endproclaim

\Proof.
We have to check that the assignments $A\mapsto H/HA^+$ and $C\mapsto\lco CH$ 
give well-defined inverse maps between the two sets. One part of the argument 
is provided by Takeuchi's Theorems 1 and 2. Let us explain the other part.

Let $C$ be any left $H$-module factor coalgebra of $H$ over which $H$ is left 
faithfully coflat, and let $A=\lco CH$. By Corollary 1.3 $C=H/HB^+$ where 
$B=H\co C$, and $H$ is right faithfully flat over $B$. By Lemma 1.6 $A=S(B)$. 
Since $S:H\to H$ is an algebra antiautomorphism, $H$ is left faithfully flat 
over $A$. By Lemma 1.4 $HA^+=HB^+$. Hence $C=H/HA^+$.

For the case of right (co)flatness we proceed similarly. Let $A$ be any right 
coideal subalgebra of $H$ over which $H$ is right faithfully flat, and let 
$C=H/HA^+$. By Corollary 1.2 $A=\lco DH$ where $D=H/A^+H$, and $H$ is left 
faithfully coflat over $D$. Since $S:H\to H$ is a coalgebra antiautomorphism, 
it follows from (1.7) that $S$ induces a coalgebra antiisomorphism $C\to D$. 
Hence $H$ is right faithfully coflat over $C$. The second equality in Lemma 
1.5 with $J=A^+H$ shows that $\lco DH=\lco CH$, i.e., $A=\lco CH$.
\endproof

The next two lemmas describe how switching of sides affects relevant 
categories of Hopf modules.

\proclaim
Lemma 1.8.
Let $C=H/HA^+$ and $D=H/A^+H$ for some right coideal subalgebra $A$ of $H$. 
There are commutative diagrams
$$ 
\diagram{
\HMC & \mapr{} & \DMH & \hskip2cm & \HCM & \mapr{} & \MHD \cr
\noalign{\smallskip}
\mapd{}{} && \mapd{}{} && \mapd{}{} && \mapd{}{} \cr
\noalign{\smallskip}
\MC & \mapr{} & \DM, && \CM & \mapr{} & \MD \cr
}
$$
where all vertical arrows are forgetful functors, while all horizontal arrows 
are category equivalences.
\endproclaim

\proclaim
Lemma 1.9.
Let $A=\lco CH$ and $B=H\co C$ for some left $H$-module factor coalgebra 
$C$ of $H$. There are commutative diagrams
$$ 
\diagram{
\AMH & \mapr{} & \HMB & \hskip2cm & \MAH & \mapr{} & \BHM \cr
\noalign{\smallskip}
\mapd{}{} && \mapd{}{} && \mapd{}{} && \mapd{}{} \cr
\noalign{\smallskip}
\AM & \mapr{} & \MB, && \MA & \mapr{} & \BM \cr
}
$$
where all vertical arrows are forgetful functors, while all horizontal arrows 
are category equivalences.
\endproclaim

\Proof.
Lemmas 1.8 and 1.9 both follow from the fact that the antipode $S$ is a Hopf 
algebra antiautomorphism of $H$ which induces an antiisomorphism, respectively, 
of coalgebras $C\to D$ and algebras $B\to A$ by Lemmas 1.4 and 1.6. For example, 
the category $\HMC$ may be identified with $\CopM_{H\opcop}$, and the latter 
is equivalent to $\DMH$ since $S$ gives an isomorphism of Hopf algebras 
$H\opcop\to H$ which induces an isomorphism of factor coalgebras $C\cop\to D$.
\endproof

\section
2. Projectivity and generation versus injectivity and cogeneration

In this section we present several conclusions which show interrelation 
between properties of right coideal subalgebras and left $H$-module factor 
coalgebras. These conclusions can be derived as special cases from more 
complex results found in literature. We will see that all of them can be 
proved by modification of an argument used in the proof of the Masuoka-Wigner 
theorem \cite{Ma-W94, Th. 2.1, (d)$\Rar$(c)}. We retain the assumptions of 
section 1.

Recall that a left (right) $C$-comodule is called \emph{cofree} if it is 
isomorphic to a direct sum of copies of $C$. It is well known that a comodule 
is injective if and only if it is coflat, if and only if it is a direct 
summand of a cofree comodule \cite{Tak77, Appendix 2}.

A $C$-comodule $U$ is a cogenerator if and only if it contains each 
indecomposable injective comodule as a direct summand. Since $C$, as a left 
(right) $C$-comodule, is a direct sum of indecomposable injectives, an 
equivalent condition is that $C$ is a direct summand of a direct sum of copies 
of $U$. The faithfully coflat comodules are precisely the injective 
cogenerators.

The Hopf algebra $H$ is always viewed as a left or right module over its 
subalgebras with respect to the action given by the multiplication in $H$, and 
as a left or right comodule over its factor coalgebras with respect to the 
coaction induced by the comultiplication in $H$.

\setitemsize(iii)
\proclaim
Proposition 2.1.
Let $A$ be a right coideal subalgebra of $H$ and $C=H/HA^+$. Then:

\item(i)
If $H$ is projective in $\AM$ then $H$ is a cogenerator in $\CM$.

\item(ii)
If $H$ is a generator in $\AM$ then $H$ is injective in $\CM$.

\item(iii)
If $H$ is projective in $\MA$ then $H$ is a cogenerator in $\MC$.

\item(iv)
If $H$ is a generator in $\MA$ then $H$ is injective in $\MC$.

\endproclaim

\Proof.
For each left $A$-module $V$ the left coaction of $C$ on $H$ gives rise to a 
coaction on $H\ot_AV$. In this way we get a functor $\Phi:\AM\to\CM$, 
forgetting part of the structure provided by Takeuchi's functor (1.2). 
Clearly, $\Phi(A)\cong H$. Let us evaluate $\Phi(H)$.

For each object $M\in\MAH$ there is a $k$-linear bijection 
$$
M\ot_AH\cong M/MA^+\ot H\eqno(2.1)
$$ 
defined by an explicit formula \cite{Tak79, p.~456}. In particular,
$$
\Phi(H)=H\ot_AH\cong H/HA^+\ot H=C\ot H.
$$
This can be also seen as follows. The image $T=H\ot_AH$ of $H\in\MAH$ under the 
functor $?\ot_AH$ is an object of $\MHH$. By the fundamental theorem on Hopf 
modules \cite{Sw, Th. 4.1.1} the right action of $H$ on $T$ yields a $k$-linear 
bijection $T\cong T\co H\ot H$. Furthermore,
$$
T\co H\cong T/TH^+\cong H/HA^+=C.
$$
The right coaction of $H$ on $T$ is the tensor product of right coactions on 
the two tensorands. We see that the left coaction of $C$ on $T$ commutes with 
both the right action and the right coaction of $H$. Hence both $T\co H$ and 
$TH^+$ are $C$-subcomodules of $T$. So too is the image of $H$ in $T$ 
under the map such that $h\mapsto h\ot1$ for $h\in H$. Therefore all 
isomorphisms in the displayed line above are left $C$-colinear. We get 
$T\co H\cong C$ in $\CM$, and $\Phi(H)\cong C\ot H$ is an isomorphism in $\CM$ 
as well. Thus $\Phi(H)$ is cofree.

If $H$ is projective in $\AM$, then $H$ is a direct summand of some free left 
$A$-module $F$. Then $\Phi(H)$ is a $\CM$-direct summand of $\Phi(F)$, whence 
so too is $C$ since $C$ is a $\CM$-direct summand of $\Phi(H)$. On the other 
hand, $\Phi(F)$ is direct sum of a family of copies of $\Phi(A)\cong H$. This 
proves (i).

If $H$ is a generator in $\AM$, then $A$ is an $\AM$-direct summand of direct 
sum of several copies of $H$. In this case $H\cong\Phi(A)$ is a $\CM$-direct 
summand of direct sum of several copies of $\Phi(H)$. Thus $H$ is a 
$\CM$-direct summand of a cofree comodule, and (ii) follows.

Put $D=H/A^+H$. Applying (i) and (ii) to the Hopf algebra $H\op$ and its right 
coideal subalgebra $A\op$, we deduce that $H$ is a cogenerator in $\DM$ under 
the hypothesis of (iii) and $H$ is injective in $\DM$ under the hypothesis of 
(iv). This amounts to the conclusions of (iii) and (iv) since $C$ is 
antiisomorphic to $D$ by (1.7).
\endproof

\proclaim
Proposition 2.2.
Let $C$ be a left $H$-module factor coalgebra of $H$ and $A=\lco CH$. Then:

\item(i)
If $H$ is injective in $\MC$ then $H$ is a generator in $\MA$.

\item(ii)
If $H$ is a cogenerator in $\MC$ then $H$ is projective in $\MA$.

\item(iii)
If $H$ is injective in $\CM$ then $H$ is a generator in $\AM$.

\item(iv)
If $H$ is a cogenerator in $\CM$ then $H$ is projective in $\AM$.

\endproclaim

\Proof.
Consider the functor $\Psi:\MC\to\MA$ that takes a right $C$-comodule $V$ to 
$$
\Psi(V)=V\sqCH\sbs V\ot H
$$
with the $A$-module structure given by right multiplications on $H$. Then 
$\Psi(C)\cong H$. Recalling $k$-linear bijections
$$
H\sqC M\cong H\ot\lco CM\eqno(2.2)
$$
constructed in \cite{Tak79, Proof of Th. 2} for Hopf modules $M\in\HCM$, we find
$$
\Psi(H)=H\sqCH\cong H\ot\lco CH=H\ot A.
$$
Under this isomorphism $h\ot a\in H\ot A$ corresponds to 
$\De(h)\cdot(1\ot a)\in H\sqCH$. Thus $\Psi(H)$ is a free $A$-module.

If $H$ is injective in the category $\MC$, then $H$ is a direct summand of 
some cofree right $C$-comodule $F$, whence the nonzero free $A$-module 
$\Psi(H)$ is a direct summand of the $A$-module $\Psi(F)$ which is a direct 
sum of copies of $H$. If $H$ is a cogenerator in $\MC$, then $C$ is an 
$\MC$-direct summand of a direct sum of copies of $H$, whence $H\cong\Psi(C)$ is 
an $\MA$-direct summand of a direct sum of copies of $\Psi(H)$. We thus deduce 
(i) and (ii). The two remaining assertions are obtained by applying (i) and 
(ii) to the Hopf algebra $H\cop$ and its left module factor coalgebra $C\cop$ 
with the aid of Lemma 1.6.
\endproof

Items (iii), (iv) of Proposition 2.2 together with Takeuchi's Theorem 1 show 
that the assumption that $H$ is left faithfully flat over its right coideal 
subalgebra $A$ imply that $H$ is a projective generator in the category $\AM$. 
This is part of Masuoka and Wigner's theorem \cite{Ma-W94, Th. 2.1}, and here 
bijectivity of the antipode is used in an essential way. Note, however, that 
(i) and (ii) of Propositions 2.1 and 2.2 do not require bijectivity of the 
antipode.

It is well known that a ring $R$ is a generator in the category $\MA$ of right 
modules over its subring $A$ if and only if $A$ is an $\MA$-direct summand of 
$R$ (see \cite{Fai, 3.27}). The dual assertion for coalgebras is also true:

\proclaim
Lemma 2.3.
Suppose that $\pi:E\to C$ is a surjective homomorphism of coalgebras. Then $E$ 
is a cogenerator in $\MC$ if and only if $\pi$ is a split epimorphism in $\MC$.
\endproclaim

\Proof.
One direction is immediately clear: if $\pi$ splits, then $C$ is an 
$\MC$-direct summand of $E$, whence $E$ is a cogenerator in $\MC$ since so is 
$C$.

Conversely, suppose that $E$ is a cogenerator in $\MC$. Then $C$ is isomorphic 
to an $\MC$-direct summand of a suitable direct sum $S$ of a family of copies 
of $E$. There are right $C$-colinear maps $f:C\to S$ and $g:S\to C$ such that 
$gf=\id$.

Denote by $\ep:C\to k$ and $\ep':E\to k$ the counits of $C$ and $E$. Recall 
that, given a right $C$-comodule $V$, the assignment $\al\mapsto\ep\circ\al$ 
gives a bijection between the set of right $C$-colinear maps $V\to C$ and the 
dual space $V^*$ of linear maps $V\to k$. A similar description of morphisms 
into $E$ is valid in $\calM^E$.

Now $S$ is a right $E$-comodule, and therefore there is a unique right 
$E$-colinear map $h:S\to E$ such that $\ep'h=\ep g$. Then $\pi h=g$ since both 
$\pi h$ and $g$ are right $C$-colinear maps $S\to C$ which correspond to the 
same linear map
$$
\ep\pi h=\ep'h=\ep g:S\to k.
$$
We have $\pi hf=gf=\id$, i.e., $hf:C\to E$ is a right $C$-colinear splitting 
of $\pi$.
\endproof

Thus the generator/cogenerator properties in Propositions 2.1 and 2.2 can be 
rephrased by saying, respectively, that the inclusion map $A\to H$ is a split 
monomorphism of $A$-modules and the canonical surjection $\pi:H\to C$ is a 
split epimorphism of $C$-comodules. In such a form both propositions are 
special cases of results due to Doi \cite{Doi92, (3.3), (3.5)}.

Still, there are some differences. In the case when $H$ is injective in $\MC$ 
it follows from \cite{Doi92, Prop. 3.3} that the left coideal subalgebra 
$B=H\co C$ is a $\BM$-direct summand of $H$. Using a \emph{total integral} 
$\ph:C\to H$, i.e., a right $C$-colinear map sending $1_C\in C$ to $1\in H$, 
Doi constructed a left $B$-linear splitting $H\to B$ of the inclusion map 
$B\to H$. It can be defined explicitly by the assignment 
$$ 
h\mapsto\sum\,h\1\,\ph\bigl(S(h\2)\cdot1_C\bigr)
=\sum\,h\1\,\ph\bigl(\pi(S(h\2))\bigr).
$$
On the other hand, part (i) of our Proposition 2.2 implies that the right 
coideal subalgebra $A=\lco CH$ is an $\MA$-direct summand of $H$. Both 
conclusions do not require bijectivity of the antipode $S$. However, they are 
not equivalent to each other unless $S$ is bijective. By the way, the 
assignment
$$
h\mapsto\sum\,S\bigl(\ph(\pi(h\1))\bigr)\,h\2
$$
defines a right $A$-linear splitting $p:H\to A$ of the inclusion map $A\to H$. 
Here the inclusion $\Img p\sbs A$ holds because 
$\sum\ph\bigl(\pi(h\1)\bigr)\ot h\2\in H\sqCH$ for each $h\in H$, while 
$H\sqCH$ is taken to the invariants $\lco CH$ by the linear map $H\ot H\to H$ 
such that $x\ot y\mapsto S(x)y$ for $x,y\in H$. Similar observations apply to 
other conclusions. Thus the technique of \cite{Doi92} provides an alternative 
approach to the proofs of Propositions 2.1, 2.2.

Let $D=H/A^+H$. The fact that a $\DM$-splitting $D\to H$ of the canonical 
surjection $H\to D$ implies that $H$ is left projective over the subalgebra 
$H\co D$ can be also extracted from the proof of \cite{Schn92, Cor. 1.6}. 
Other closely related statements are found in \cite{Scha-Schn05, Cor. 2.7, 
Remark 4.3}.

\section
3. Tensoring operations

In this section we recall tensoring operations on the categories of comodules, 
modules, and Hopf modules. They provide a useful tool for later discussions. 
Here we present results only in the form which does not require bijectivity of 
the antipode of $H$. Let $A$ be a right coideal subalgebra and $C$ a left 
$H$-module factor coalgebra of the Hopf algebra $H$.

Except for Propositions 3.1 and 3.5 we may assume even more generally that $A$ 
is any right $H$-comodule algebra and $C$ any left $H$-module coalgebra.

For objects $U\in\calM^H$ and $V\in\MC$ the tensor product $U\ot V$ over the base 
field is a right $(H\ot C)$-comodule. Since the $H$-module structure map 
$H\ot C\to C$ is a homomorphism of coalgebras, it equips $U\ot V$ with a 
right $C$-comodule structure
$$
u\ot v\mapsto\sum\,(u\0\ot v\0)\ot u\1v\1,\qquad u\in U,\ v\in V.\eqno(3.1)
$$
If $M\in\MAH$, then also $U\ot M\in\MAH$ with respect to the given action of 
$A$ on $M$ and the tensor product of $H$-comodule structures on $U$ and $M$.

In this way $\MC$ and $\MAH$ become left module categories over the monoidal 
category $\calM^H$ of right $H$-comodules. If $\dim U<\infty$, then there is a 
right $H$-comodule structure on the dual space $U^*$ which makes $U^*$ the 
left dual of $U$ in the category $\calM^H$. It is immediate then that $U^*\ot\,?$ 
and $U\ot\,?$, regarded as endofunctors of either $\MC$ or $\MAH$, are a pair 
of adjoint functors, i.e., there are natural bijections between the sets of 
morphisms in the respective categories
$$
\openup1\jot
\eqalignno{
\MC(V,\,U\ot W)&\cong\MC(U^*\ot V,\,W),&(3.2)\cr
\MAH(M,\,U\ot N)&\cong\MAH(U^*\ot M,\,N)&(3.3)
}
$$
for objects $V,W\in\MC$ and $M,N\in\MAH$.

\proclaim
Proposition 3.1.
Suppose that $A\sbs\lco CH$. Takeuchi's functors $\Phi:\MAH\to\MC$ and 
$\Psi:\MC\to\MAH$ commute with the tensoring functors $U\ot\,?$, i.e., 
there are natural isomorphisms
$$
\Phi(U\ot M)\cong U\ot\Phi(M),\qquad
\Psi(U\ot V)\cong U\ot\Psi(V)
$$
for objects $U\in\calM^H,$ $\,V\in\MC,$ and $M\in\MAH$. In particular,
$\,\Psi(U)\cong U\ot\lco CH\,$ where $U\in\calM^H$ is regarded as a right 
$C$-comodule.
\endproclaim

\Proof.
We have $k$-linear bijections
$$
\Phi(U\ot M)=(U\ot M)/(U\ot MA^+)\cong U\ot M/MA^+=U\ot\Phi(M)
$$
which are $C$-colinear since the $C$-comodule structures on both the left and 
right hand vector spaces arise from the $H$-comodule structure on $U\ot M$.

If $\dim U<\infty$, then the functors $V\mapsto\Psi(U\ot V)$ and 
$V\mapsto U\ot\Psi(V)$ are right adjoints, respectively, of the functors 
$M\mapsto U^*\ot\Phi(M)$ and $M\mapsto\Phi(U^*\ot M)$. The last two functors 
are isomorphic, as we have just seen. Hence so are their right adjoints. This 
shows that $\Psi(U\ot V)\cong U\ot\Psi(V)$ when $\dim U<\infty$. Since these 
isomorphisms are natural in $U$, it follows from the local finite-dimensionality 
of comodules that the isomorphisms extend to arbitrary right $H$-comodules $U$. 
One has just to note that the functor $\Psi={?\sqCH}$ preserves directed 
colimits \cite{Br-W, 10.5}, and so do the tensor products over the base field.

Now take $V=k1_C$. This is the 1-dimensional $C$-comodule corresponding to the 
distinguished grouplike $1_C\in C$. Then $U\ot V\cong U$ in $\MC$ and 
$\Psi(V)\cong\lco CH$ in $\MAH$, whence the final assertion.
\endproof

Since $\Psi(V)=V\sqCH$, Proposition 3.1 gives an isomorphism of $H$-comodules
$$
(U\ot V)\sqCH\cong U\ot(V\sqCH)\eqno(3.4)
$$
for $U\in\calM^H$ and $V\in\MC$. In the case when $H$ is commutative and $G$ is 
the affine group scheme represented by $H$ the functor $?{}\sqCH$ is the 
induction functor from the representation category of a subgroup of $G$ to 
the representation category of $G$. In this case the above isomorphism is a 
classical fact known as the Tensor Identity \cite{Ja, Prop. 3.6}. The quantum 
group variants of the tensor identity are discussed by Parshall and Wang 
\cite{Par-W91, Th. 2.7.1}.

For later applications in sections 8 and 9 we will need a more general version 
of formula (3.4) which will be presented in Proposition 3.4.

The category $\CM$ of left $C$-comodules is likewise a left module category 
over the monoidal category $\HcoM$ of left $H$-comodules. Transformation by 
the antipode $S$ gives rise to a functor $\calM^H\to\HcoM$. It allows us to view 
$\CM$ as a right module category over the monoidal category $\calM^H$. Explicitly, 
for $U\in\calM^H$ and $W\in\CM$ the left coaction of $C$ on the vector space 
$W\ot U$ is defined by the rule
$$
w\ot u\mapsto\sum S(u\1)w_{(-1)}\ot(w\0\ot u\0),\qquad u\in U,\ w\in W.\eqno(3.5)
$$
Denote by $U\triv$ the vector space $U$ equipped with the trivial coaction of 
$H$ such that $u\mapsto u\ot1$ for all $u\in U$. Tensoring with $U\triv$ results 
in the vector space $U\ot V$ with the given coaction of $C$ on $V$ not affected 
by $U$.

\proclaim
Lemma 3.2.
Given $U\in\calM^H,$ $\,V\in\HMC,$ and $\,W\in\HCM,$ we have
$$
U\ot V\cong U\triv\ot V\ \ \text{in $\,\MC$},\qquad 
W\ot U\cong W\ot U\triv\ \ \text{in $\,\CM$}.
$$
\endproclaim

\Proof.
The linear transformation of the vector space $U\ot V$ such that
$$
u\ot v\mapsto\sum u\0\ot u\1v\quad\text{for $u\in U$ and $v\in V$}
$$
gives an isomorphism of $C$-comodules $U\ot V\to U\triv\ot V$. The inverse 
isomorphism sends $u\ot v$ to $\,\sum u\0\ot S(u\1)v$. Similarly, the assignment
$$
w\ot u\mapsto\sum u\1w\ot u\0,\qquad u\in U,\ w\in W,
$$
defines an isomorphism of $C$-comodules $W\ot U\triv\to W\ot U$.
\endproof

\proclaim
Corollary 3.3.
The following classes of right and left $C$-comodules are stable under the 
functors, respectively, $U\ot\,?$ and $?{}\ot U$ for all $U\in\calM^H$: cofree 
comodules, injective comodules, comodules of finite injective dimension.
\endproclaim

\Proof.
Since $C\in\HMC$, the $C$-comodule $U\ot C\cong U\triv\ot C$ is cofree. It 
follows that $U\ot V$ is cofree whenever so is $V\in\MC$. Since the functor 
$U\ot\,?$ is additive, it preserves injective comodules which are precisely 
direct summands of cofree ones. Since this functor is exact, the 
injective dimension of $U\ot V$ for any $V\in\MC$ cannot exceed the injective 
dimension of $V$. All this applies to $\CM$ as well.
\endproof

\proclaim
Proposition 3.4.
There are natural isomorphisms of vector spaces
$$
(U\ot V)\sqC W\cong V\sqC(W\ot U)\eqno(3.6)
$$
for $\,U\in\calM^H,$ $V\in\MC,$ and $W\in\CM$. Moreover,
$$
(U\ot V)\sqC W\cong V\sqC(W\ot U)\cong U\ot(V\sqC W)\eqno(3.7)
$$
if either $\,V\in\HMC$ or $\,W\in\HCM$.
\endproclaim

\Proof.
Let $\rho:V\to V\ot C$ and $\la:W\to C\ot W$ be the comodule structure maps. 
Denote by $\rho_U$ and $\la_U$ the comodule structure maps on the vector 
spaces $U\ot V$ and $W\ot U$ defined, respectively, by formulas (3.1) and (3.5). 
Then $(U\ot V)\sqC W$ and $V\sqC(W\ot U)$ coincide with the equalizers 
of the pairs of linear maps
$$
\eqalign{
&\rho_U\ot\id_W,\ \id_U\ot\id_V\ot\la:\ U\ot V\ot W\to U\ot V\ot C\ot W
\qquad\text{and}\cr
&\rho\ot\id_W\ot\id_U,\ \id_V\ot\la_U:\ V\ot W\ot U\to V\ot C\ot W\ot U,
}
$$
respectively. Let $\al:V\ot W\ot U\to U\ot V\ot W$ be the $k$-linear bijection 
obtained by permutation of tensorands. Define a $k$-linear map
$$
\be:V\ot C\ot W\ot U\to U\ot V\ot C\ot W
$$
by the rule $\be(v\ot c\ot w\ot u)=\sum u\0\ot v\ot u\1c\ot w\,$ for $u\in U$, 
$v\in V$, $w\in W$, and $c\in C$. As one checks straightforwardly, 
$$
\eqalign{
(\rho_U\ot\id_W)\circ\al&=\be\circ(\rho\ot\id_W\ot\id_U),\cr
(\id_U\ot\id_V\ot\la)\circ\al&=\be\circ(\id_V\ot\la_U).
}
$$
Since $\al$ and $\be$ are both invertible, it follows that $\al$ maps 
$V\sqC(W\ot U)$ bijectively onto $(U\ot V)\sqC W$, yielding (3.6).
If either $\,V\in\HMC$ or $\,W\in\HCM$, then we may replace $U$ by $U\triv$, 
making use of Lemma 3.2. This space is not involved in the coactions of $C$, 
so we get (3.7).
\endproof

Consider now another pair of functors (1.2). For modules $U\in\HM$ and 
$V\in\AM$ the tensor product $V\ot U$ is a left $A\ot H$-module. It is then 
a left $A$-module with $A$ acting via the algebra homomorphism $A\to A\ot H$ 
given by the $H$-comodule structure on $A$. In Sweedler's notation
$$
a\cdot(v\ot u)=\sum a\0v\ot a\1u,\qquad a\in A,\ v\in V,\ u\in U.\eqno(3.8)
$$
If $M\in\HCM$, then also $M\ot U\in\HCM$ with respect to the given coaction of 
$C$ on $M$ and the tensor product of $H$-module structures on $M$ and $U$. 

In this way $\AM$ and $\HCM$ are right module categories over the monoidal 
category $\HM$ of left $H$-modules. We will view the category $\MA$ of right 
$A$-modules as a left module category over $\HM$, defining for each $U\in\HM$ 
and $W\in\MA$ a right action of $A$ on the vector space $U\ot W$ by the rule
$$
(u\ot w)\cdot a=\sum S(a\1)u\ot wa\0,\qquad a\in A,\ u\in U,\ w\in W.\eqno(3.9)
$$

\proclaim
Proposition 3.5.
Suppose that $A\sbs\lco CH$. Takeuchi's functors $\Phi:\AM\to\HCM$ and 
$\Psi:\HCM\to\AM$ commute with the tensoring functors $?\ot U$, i.e., 
there are natural isomorphisms
$$
\Phi(V\ot U)\cong\Phi(V)\ot U,\qquad
\Psi(M\ot U)\cong\Psi(M)\ot U
$$
for objects $U\in\HM,$ $\,V\in\AM,$ and $M\in\HCM$. In particular,
$\,\Phi(U)\cong H/HA^+\ot U\,$ where $U\in\HM$ is regarded as a left 
$A$-module.
\endproclaim

\Proof.
The second isomorphism $\lco C(M\ot U)\cong\lco CM\ot U$ is clear since $U$ is 
not involved in the coaction of $C$. The first isomorphism is obtained by 
means of the $H$-linear map
$$
\openup1\jot
\displaylines{
H\ot_A(V\ot U)\to(H\ot_AV)\ot U,\cr
h\ot(v\ot u)\mapsto\sum\,(h\1\ot v)\ot h_2u,\qquad h\in H,\ v\in V,\ u\in U.
}
$$
Since $H\in\MAH$, this is a special case of Proposition 3.8. Finally, 
for $V=A/A^+$ we have $V\ot U\cong U$ and $\Phi(V)\cong H/HA^+$.
\endproof

In the next lemma we denote by $U\triv$ the vector space $U$ equipped with the 
trivial action of $H$ such that $hu=\ep(h)u$ for all $u\in U$ and $h\in H$.

\proclaim
Lemma 3.6.
Given $U\in\HM,$ $\,V\in\AMH,$ and $\,W\in\MAH,$ we have
$$
V\ot U\cong V\ot U\triv\ \ \text{in $\,\AM$},\qquad 
U\ot W\cong U\triv\ot W\ \ \text{in $\,\MA$}.
$$
\endproclaim

\Proof.
Invertible $A$-linear maps $V\ot U\triv\to V\ot U$ and $U\ot W\to U\triv\ot W$ 
are defined, respectively, by the assignments
$$
v\ot u\mapsto\sum v\0\ot v\1u,\qquad
u\ot w\mapsto\sum w\1u\ot w\0
$$
for $u\in U$, $v\in V$, and $w\in W$.
\endproof

\proclaim
Corollary 3.7.
The following classes of right and left $A$-modules are stable under the 
functors, respectively, $U\ot\,?$ and $?{}\ot U$ for all $U\in\HM$: free 
modules, projective modules, flat modules, modules of finite projective 
dimension, modules of finite flat dimension.
\endproclaim

\Proof.
Any free $A$-module $F$ admits an $H$-comodule structure which makes $F$ 
an object of either $\AMH$ or $\MAH$. Hence $U\ot F\cong U\triv\ot F$ is a 
free right $A$-module whenever so is $F$. Similarly, the functor $?{}\ot U$ 
preserves freeness of left $A$-modules. Projective modules are direct summands 
of free modules, while flat modules are directed colimits of free modules by 
Lazard's theorem. Since tensoring functors are additive and commute with 
directed colimits, they preserve those properties of $A$-modules. For the 
remaining two classes of $A$-modules the conclusion follows by exactness of 
tensoring functors.
\endproof

\proclaim
Proposition 3.8.
There are natural isomorphisms of vector spaces
$$
(U\ot W)\ot_AV\cong W\ot_A(V\ot U)\eqno(3.10)
$$
for $\,U\in\HM,$ $\,V\in\AM,$ and $\,W\in\MA$. Moreover,
$$
(U\ot W)\ot_AV\cong W\ot_A(V\ot U)\cong(W\ot_AV)\ot U\eqno(3.11)
$$
if either $\,V\in\AMH$ or $\,W\in\MAH$.
\endproclaim

\Proof.
For all $u\in U$, $v\in V$, $w\in W$, and $a\in A$ we have
$$
\eqalign{
w\ot(av\ot u)&=\sum w\ot\bigl(a\0v\ot a\1S(a\2)u\bigr)
=\sum wa\0\ot\bigl(v\ot S(a\1)u\bigr),\cr
(u\ot wa)\ot v&=\sum\,\bigl(S(a\1)a\2u\ot wa\0\bigr)\ot v
=\sum\,(a\1u\ot w)\ot a\0v,
}
$$
respectively, in $W\ot_A(V\ot U)$ and $(U\ot W)\ot_AV$. It follows that there 
are mutually inverse well-defined linear maps between the two vector spaces in 
(3.10) under which $(u\ot w)\ot v$ corresponds to $w\ot(v\ot u)$. If 
$\,V\in\AMH$ or $\,W\in\MAH$, then Lemma 3.7 allows us to replace $U$ by 
$U\triv$, and (3.11) reduces to the standard associativity of tensor products.
\endproof

\section
4. Characterizing flatness over coideal subalgebras

In the spirit of Propositions 2.1, 2.2 it will be shown that the flatness 
property of $H$ over $A$ corresponds to the property that each finite-dimensional 
$C$-comodule admits a colinear embedding in some $H$-comodule. While the 
former property is a weakening of projectivity, the latter may be viewed as a 
weakened cogenerator property of $H$: the elements of any finite-dimensional 
$C$-comodule $V$ are separated by $C$-colinear maps $V\to H$. Precise statements 
given in Propositions 4.2 and 4.3 lead to a bijection of Theorem 4.6.

We assume that $A$ is a right coideal subalgebra of $H$ and $C$ a left 
$H$-module factor coalgebra of $H$. 

\setitemsize(a)
\proclaim
Lemma 4.1.
Suppose that $C=H/HA^+$. Consider the following properties:

\item(a)
$H$ is left flat over $A$,

\item(b)
Takeuchi's functor $\,\Phi:\MAH\to\MC\,$ is exact,

\item(c)
Takeuchi's functor $\,\Psi:\MC\to\MAH\,$ is fully faithful,

\item(d)
each right $C$-comodule is a $C$-homomorphic image of some right $H$-comodule.

Then $\,\mathrm{(a)}\Lrar\mathrm{(b)}\Rar\mathrm{(c)}\Rar\mathrm{(d)}\,$ where 
only $\,\mathrm{(b)}\Rar\mathrm{(a)}\,$ requires bijectivity of the antipode. 
Also, $H$ is left faithfully flat over $A$ if and only if\/ $\Phi$ is 
faithfully exact.
\endproclaim

\Proof.
Recall that $\Phi(M)=M/MA^+$. It follows from (2.1) that $\Phi$ is exact 
if and only if the functor $?\ot_AH$ is exact on the category $\MAH$. Hence 
$\mathrm{(a)}\Rar\mathrm{(b)}$. Also, $\Phi(M)=0$ if and only if $M\ot_AH=0$. 
Hence faithful flatness of $H$ over $A$ implies that $\,\Ker\Phi=0$, i.e., 
$\Phi$ is faithfully exact. The converse is proved in \cite{Sk10, Lemma 1.2}. 
The argument given there is based on the observation that for each right 
$A$-module $V$ one makes $V\ot H$ into an object of the category $\MAH$ in 
such a way that there is a natural isomorphism of vector spaces 
$\,V\ot_A\!H\cong\Phi(V\ot H)$.

The implication $\mathrm{(b)}\Rar\mathrm{(c)}$ is proved in \cite{Sk10, Lemma 1.1}. Full 
faithfulness of $\Psi$ means that the adjunction $\eta_V:\Phi\Psi(V)\to V$ is 
an isomorphism in $\MC$ for each right $C$-comodule $V$. Setting 
$W=\Psi(V)=V\sqCH$, we get a right $H$-comodule $W$ and a $C$-colinear map 
$W\to V$ which factors as the composite of two surjections 
$$
W\to W/WA^+=\Phi\Psi(V)\to V.
$$
Thus $\mathrm{(c)}\Rar\mathrm{(d)}$.
\endproof

\setitemsize(ii)
\proclaim
Proposition 4.2.
Suppose that $C=H/HA^+$ and $H$ is left (respectively, right) flat over $A$. 
Then each finite-dimensional left (respectively, right) $C$-comodule embeds as 
a $C$-subcomodule in some left (respectively, right) $H$-comodule.  
\endproclaim

\Proof.
Suppose that $H$ is left flat over $A$, and let $V$ be a finite-dimensional left 
$C$-comodule. The categories of finite-dimensional left and right $C$-comodules 
are antiequivalent under the functor which takes $V$ to the dual vector space 
$V^*$. By Lemma 4.1 there exist a right $H$-comodule $U$ and a $C$-colinear 
surjection $U\to V^*$. Since each comodule is the sum of its finite-dimensional 
subcomodules, we can find a finite-dimensional comodule $U$ satisfying the 
previous condition. Then $V$ embeds in the left $H$-comodule $U^*$.

Set $D=H/A^+H$. If $H$ is right flat over $A$, then $H\op$ is left flat over 
$A\op$. As we have just proved, in this case each finite-dimensional left 
$\DM$-comodule embeds in some left $H$-comodule. By (1.7) this translates into 
the desired property of finite-dimensional right $C$-comodules.
\endproof

\setitemsize(iii)
\proclaim
Proposition 4.3.
Suppose that  $A=\lco CH$ and each finite-dimensional right (respectively, 
left) $C$-comodule embeds as a $C$-subcomodule in some right (respectively, 
left) $H$-comodule. Then:

\item(i)
$H$ is right (respectively, left) flat over $A$.

\item(ii)
For each Hopf module $M$ in $\HCM$ (respectively, in $\HMC$) we have $M=HJ$ 
where $J=\lco CM$ (respectively, $J=M\co C$).

\item(iii)
$C=H/HA^+$.

\endproclaim

\Proof.
Consider the case when the hypothesis is satisfied for right $C$-comodules. 
The other case will follow by switching over to the Hopf algebra $H\cop$ and 
making use of Lemmas 1.4, 1.6.

\smallskip
(i) We will use Takeuchi's functor $\Psi:\MC\to\MAH$, $\,V\mapsto V\sqCH$. 
Let us say that a right $C$-comodule is \emph{$H$-liftable} if it belongs 
to the image of the canonical functor $\calM^H\to\MC$, and a right 
$C$-comodule is \emph{locally $H$-liftable} if each its finite subset is 
contained in an $H$-liftable subcomodule.

If $U\in\MC$ is $H$-liftable then $\Psi(U)\cong U\ot A$ is a free right 
$A$-module by Proposition 3.1. If $W\in\MC$ is locally $H$-liftable then the 
set $\calL$ of all its $H$-liftable finite-dimensional subcomodules is directed 
by inclusion and the union of the subcomodules in this set gives the whole $W$. 
Since $\Psi$ preserves directed colimits (see \cite{Br-W, 10.5}), we get
$$
\Psi(W)\cong\limdir_{U\in\calL}\Psi(U),
$$
i.e., $\Psi(W)$ is a directed colimit of free $A$-modules, which implies that 
this $A$-module is flat.

Suppose that $C$ embeds in some locally $H$-liftable right $C$-comodule $W$ as 
a subcomodule. Then $C$ is an $\MC$-direct summand of $W$ since $C$ is injective 
in $\MC$, and it follows that $\Psi(C)$ is a direct summand of $\Psi(W)$. In 
this case $H\cong\Psi(C)$ is a flat right $A$-module since so is $\Psi(W)$.

So it remains to prove that each right $C$-comodule $V$ can be embedded in a 
locally $H$-liftable right $C$-comodule $W$. Denote by $\calF(V)$ the set of 
all finite-dimensional subcomodules of $V$, and for each $U\in\calF(V)$ choose 
any right $H$-comodule $\widetilde U$ which contains $U$ as a $C$-subcomodule. 
The $C$-comodule $\bigoplus_{U\in\calF(V)}U$ is a subcomodule of 
$\bigoplus_{U\in\calF(V)}\widetilde U$. There is an epimorphism 
$p:\bigoplus_{U\in\calF(V)}U\to V$ in $\MC$ such that $p|_U$ is the 
inclusion map $U\to V$ for each $U\in\calF(V)$, and it follows that $V$ embeds 
in the right $C$-comodule 
$$
\textstyle V_1=\bigl(\bigoplus_{U\in\calF(V)}\widetilde U\bigr)/\Ker p.
$$
Since each direct summand $\widetilde U$ has zero intersection with $\Ker p$, 
it embeds in $V_1$ as a subcomodule. Thus for each $U\in\calF(V)$ there exists 
an $H$-liftable subcomodule $\widetilde U\sbs V_1$ containing $U$.

Starting with $V_0=V$ and iterating the preceding construction we get a 
sequence of right $C$-comodules $V_0,V_1,V_2,\ldots$ and monomorphisms 
$V_0\to V_1\to V_2\to\ldots$ in $\MC$ with the property that for each $i$ and 
each $U\in\calF(V_i)$ there exists an $H$-liftable subcomodule 
$\widetilde U\sbs V_{i+1}$ containing the image of $U$ in $V_{i+1}$. Now put 
$W=\limdir V_i$. Each $V_i$ embeds in $W$ as a subcomodule. Identifying $V_i$ 
with its image in $W$ we will have $W=\bigcup V_i$. If $U\in\calF(W)$, then 
$U\in\calF(V_i)$ for some $i$, whence there exists an $H$-liftable subcomodule 
$\widetilde U\sbs V_{i+1}\sbs W$ containing $U$. Thus $W$ is a locally 
$H$-liftable right $C$-comodule in which $V$ is embedded.

\smallskip
(ii) By duality each finite-dimensional left $C$-comodule is a homomorphic image 
of a comodule in the image of the canonical functor $\HcoM\to\CM$. Since that 
functor has a right adjoint $H\sqC{}?$, it follows that for a left 
$C$-comodule $V$ the adjunction
$$
\eta_V:H\sqC V\to V
$$
is surjective whenever $\,\dim V<\infty$. By the local finite-dimensionality 
of comodules $\eta_V$ is surjective for any $V$.

If $M\in\HCM$, then $\tilde M=H\sqCM$ is an object of $\HHM$. By \cite{Sw, Th. 
4.1.1} the subspace of coaction invariants $\lco H\tilde M$ freely generates 
$\tilde M$ as an $H$-module. Since the map $\eta_M:\tilde M\to M$ is 
$C$-colinear, we have $\eta_M(\lco H\tilde M)\sbs\lco CM=J$, and since 
$\eta_M$ is also $H$-linear, its surjectivity yields
$$
M=\eta_M(H\cdot\lco H\tilde M)=HJ.
$$

\smallskip
(iii) We have $C=H/I$ where $I$ is the kernel of the canonical surjection 
$H\to C$. So $I$ is an $\HCM$-subobject of $H$, and (ii) yields 
$I=HJ$ where
$$
J=\lco CI=I\cap\lco CH=I\cap A.
$$
Finally, $I\cap A=A^+$ by Lemma 1.1.
\endproof

Note that the unparenthesized variants of Propositions 4.2 and 4.3 do not 
require bijectivity of the antipode. Further on we will use the notion of 
a dominion subalgebra explained in the introduction.

\proclaim
Lemma 4.4.
A right coideal subalgebra $A$ of $H$ is a dominion subalgebra if and only if 
$A=\lco CH$ for some left $H$-module factor coalgebra $C$ of $H,$ and in this 
case the equality holds with $C=H/HA^+$.
\endproclaim

\Proof.
Let $\pi:H\to C$ be the canonical surjection. If $A=\lco CH$, then 
$HA^+\sbs\Ker\pi$ by Lemma 1.1. This entails $A\sbs\lco{C'}H\sbs\lco CH$, and 
therefore $\lco{C'}H=A$, where $C'=H/HA^+$. Thus it suffices to consider the 
factor coalgebra $C=H/HA^+$.

Then there is a linear bijection $H\ot_AH\cong C\ot H$ under which 
$x\ot y\in H\ot_AH$ corresponds to $\pi(x\1)\ot x\2y\in C\ot H$ 
(see \cite{Tak79, p.~456}). Therefore $h\in H$ satisfies $h\ot1=1\ot h$ in 
$H\ot_AH$ if and only if $\pi(h\1)\ot h\2=1_C\ot h$ in $C\ot H$. The latter 
equality defines the subalgebra of invariants $\,\lco CH$. Thus this 
subalgebra is precisely the dominion of $A$ in $H$.
\endproof

\proclaim
Corollary 4.5.
Let $A$ be a right coideal subalgebra of $H$. If $H$ is a generator in either 
$\AM$ or $\MA$, then $A=\lco CH$ where $C=H/HA^+$.
\endproclaim

\Proof.
There exists a left or right $A$-linear map $f:H\to A$ such that $f|_A=\id$. 
If some element $h\in H$ satisfies $h\ot1=1\ot h$ in $H\ot_AH$, then, applying 
either $\id\ot f$ or $f\ot\id$ to both sides of that equality, we find 
$h=f(h)1\in A$. In other words, $A$ is a dominion subalgebra.
\endproof

\proclaim
Theorem 4.6.
Takeuchi's correspondence gives a bijection between the set of right coideal 
dominion subalgebras of $H$ over which $H$ is left (respectively, right) flat 
and the set of left $H$-module factor coalgebras $C$ of $H$ with the property 
that each finite-dimensional left (respectively, right) $C$-comodule embeds 
as a $C$-subcomodule in a left (respectively, right) $H$-comodule.
\endproclaim

\Proof.
Denote by $\Si_1$ and $\Si_2$, respectively, the sets of right coideal 
dominion subalgebras and left $H$-module factor coalgebras of $H$ 
satisfying the properties in the statement of this theorem. Suppose that 
$A\in\Si_1$ and $C=H/HA^+$. Then $C\in\Si_2$ by Proposition 4.2 and 
$A=\lco CH$ by Lemma 4.4.

Conversely, suppose that $C\in\Si_2$ and $A=\lco CH$. Then $A\in\Si_1$ and 
$C=H/HA^+$ by Proposition 4.3.
\endproof

Dualizing, we define \emph{dominion factor coalgebras} of $H$ by the 
condition that the canonical surjection $H\to C$ is the coequalizer of the two 
maps $H\sqCH\rightrightarrows H$ that are the restrictions of the maps 
$\ep\ot\id$ and $\id\ot\ep:H\ot H\to H$.

A left $H$-module factor coalgebra $C$ of $H$ is dominion if and only if 
$C=H/HA^+$ for some right coideal subalgebra $A$ of $H$, if and only if the 
preceding equality holds with $A=\lco CH$. In the language of 
\cite{Scha-Schn05} this can be rephrased by saying that the algebra $H$ is a 
left $C$-Galois extension of its subalgebra $\lco CH$. The dual of Theorem 
4.6 follows immediately from Propositions 2.1, 2.2 and Corollary 4.5:

\proclaim
Theorem 4.7.
Takeuchi's correspondence gives a bijection between the set of right coideal 
subalgebras $A$ of $H$ such that $H$ is a generator in the category of left 
(respectively, right) $A$-modules and the set of left $H$-module dominion 
factor coalgebras of $H$ over which $H$ is left (respectively, right) coflat.  
\endproclaim

It seems that no examples are known of $H$-module factor coalgebras over which 
$H$ is coflat but not faithfully coflat. So Theorem 4.7 currently has no 
applications in contrast to Theorem 4.6.

\section
5. Several extras

The equivalence $\mathrm{(a)}\Lrar\mathrm{(b)}$ of Lemma 4.1 leads to 
interesting conclusions which are presented in this section along with their 
dual versions. However, for later application in this paper we will need only 
Corollaries 5.2, 5.7 and Lemma 5.4.

\proclaim
Corollary 5.1.
Suppose that $A\sbs H\sbs K$ where $K$ is a Hopf algebra, $H$ a Hopf 
subalgebra with bijective antipode, and $A$ a right coideal subalgebra. 
If $K$ is left (faithfully) flat over $A,$ then so too is $H$.
\endproclaim

\Proof.
By $\mathrm{(a)}\Rar\mathrm{(b)}$ of Lemma 4.1 Takeuchi's functor 
$\,\Phi:\MAK\to\calM^{K/KA^+}\,$ is (faithfully) exact. The category $\MAH$ is 
a full subcategory of $\MAK$ and the inclusion functor $\MAH\to\MAK$ is exact.  
Hence the restriction of $\Phi$ to $\MAH$ is (faithfully) exact as well, and 
the desired conclusion follows from $\mathrm{(b)}\Rar\mathrm{(a)}$ of Lemma 4.1.
\endproof

\proclaim
Corollary 5.2.
Suppose that $H$ is left flat over its right coideal subalgebra $A$. Let 
$C=H/HA^+,$ and let $A_1$ be any right coideal subalgebra of $H$ such that 
$$
A\sbs A_1\sbs\lco CH.
$$
Then $A_1^+=A_1A^+$ and $H$ is left flat over $A_1$.
\endproclaim

\Proof.
We may regard $A_1$ and $H$ as objects of the category $\MAH$. Since Takeuchi's 
functor $\,\MAH\to\MC\,$ is exact by Lemma 4.1, it takes the inclusion 
$A_1\to H$ to a monomorphism in $\MC$. In other words, the induced map 
$A_1/A_1A^+\to H/HA^+$ is injective. But $A_1^+\sbs HA^+$ by Lemma 1.1. Hence
$$
A_1^+\sbs A_1\cap HA^+=A_1A^+.
$$
The opposite inclusion $A_1A^+\sbs A_1^+$ is obvious. It follows that 
$MA_1^+=MA^+$ for every right $A_1$-module. In particular, $H/HA_1^+=H/HA^+=C$. 
Takeuchi's functor $\calM_{A_1}^H\to\MC$ factors as the composite of 
the functor $\calM_{A_1}^H\to\MAH$ obtained by viewing $A_1$-modules as 
$A$-modules and Takeuchi's functor $\MAH\to\MC$. Since the latter two functors 
are both exact, so is their composite, whence $H$ is left flat over $A_1$ by 
Lemma 4.1.
\endproof

\setitemsize(iii)
\proclaim
Corollary 5.3.
Suppose that $H$ is left flat over its subbialgebra $A$. Let\/ $C=H/HA^+\!$. 
Then:

\item(i)
$\lco CH$ is the Hopf subalgebra of $H$ generated by $A$.

\item(ii)
$\lco CH=S(\lco CH)$ and $H$ is both left and right faithfully flat over $\lco CH$.

\item(iii)
$H$ is left faithfully flat over $A$ if and only if $A$ is a Hopf subalgebra.

\endproclaim

\Proof.
Denote by $A_1$ the subalgebra of $H$ generated by $\bigcup_{n\ge0}S^n(A)$. 
Since each $S^n(A)$ is a subbialgebra of $H$, so too is $A_1$. Also, 
$S(A_1)\sbs A_1$. It is clear therefore that $A_1$ is the smallest Hopf 
subalgebra of $H$ containing $A$. For each $n$ we have
$$
HS^n(A)^+=HS^{n+1}(A)^+
$$
by Lemma 1.4. Hence $S^n(A)^+\sbs HA^+$ by induction, and $S^n(A)\sbs\lco CH$ 
by Lemma 1.1. It follows that $A_1\sbs\lco CH$. By Corollary 5.2 $H$ is left 
flat over $A_1$. But $A_1$ is a simple object of the category $\calM_{A_1}^H$ 
since $A_1$ is a Hopf algebra. Hence $H$ is left faithfully flat over $A_1$ by 
the Masuoka-Wigner theorem \cite{Ma-W94}, and $A_1=\lco CH$ by Takeuchi's 
Theorem 1. Now $S^{-1}(A_1)$ is a subbialgebra of $H$. Since
$$ 
HS^{-1}(A_1)^+=HA_1^+=HA^+
$$
by Lemma 1.4, we get $S^{-1}(A_1)\sbs\lco CH$ by Lemma 1.1. This yields 
$S(A_1)=A_1$. Since $S$ is an antiautomorphism of $H$, faithful flatness 
over $A_1$ holds also on the right. Both (i) and (ii) have been verified. 
If $H$ is left faithfully flat over $A$, then $A=\lco CH$ by Takeuchi's 
Theorem 1. So (iii) follows from (i) and (ii).
\endproof

Note that assertion (iii) of Corollary 5.3 has long been known. 
The fact that faithful flatness of a Hopf algebra over its subbialgebra $A$ 
implies that $A$ is necessarily a Hopf subalgebra was observed by Takeuchi 
\cite{Tak94, Cor. 1.5} (see also \cite{Scha00, Cor.  2.4}). The opposite 
direction in (iii) is provided by the Masuoka-Wigner theorem \cite{Ma-W94}.

The next result is a dual version of Lemma 4.1:

\proclaim
Lemma 5.4.
Let $C$ be a left $H$-module factor coalgebra of $H,$ and let $A=\lco CH$. 
Consider the following properties:

\item(a)
$H$ is right coflat over $C$,

\item(b)
Takeuchi's functor $\,\Psi:\HCM\to\AM\,$ is exact,

\item(c)
Takeuchi's functor $\,\Phi:\AM\to\HCM\,$ is fully faithful,

\item(d)
each left $A$-module embeds as an $A$-submodule in some left $H$-module.

Then $\,\mathrm{(a)}\Lrar\mathrm{(b)}\Rar\mathrm{(c)}\Rar\mathrm{(d)}\,$ where 
only $\,\mathrm{(b)}\Rar\mathrm{(a)}\,$ requires bijectivity of the antipode. 
Also, $H$ is right faithfully coflat over $C$ if and only if\/ $\Psi$ is 
faithfully exact.
\endproclaim

\Proof.
Here $\Psi(M)=\lco CM$, and it follows from (2.2) that $\Psi$ is exact if and 
only if the functor $H\sqC{}?$ is exact on the category $\HCM$. Hence 
$\mathrm{(a)}\Rar\mathrm{(b)}$.

We will prove next that $\mathrm{(b)}\Rar\mathrm{(a)}$. Let $V$ be a left 
$C$-comodule. Then $H\ot V$ is an object of $\HCM$ with respect to the action 
of $H$ by left multiplications on the first tensorand and the left $C$-comodule 
structure $H\ot V\to C\ot H\ot V$ such that
$$
h\ot v\mapsto\sum\,h\1v_{(-1)}\ot h\2\ot v\0,\qquad h\in H,\ v\in V.
$$
We have
$$
\Psi(H\ot V)=\lco C(H\ot V)\cong k1_C\sqC(H\ot V).
$$
On the other hand, we have defined in section 3 tensoring operations by means 
of right $H$-comodules. In particular, $V\ot H$ is a left $C$-comodule with 
respect to the coaction of $C$ defined by the rule
$$
v\ot h\mapsto\sum\,S(h\2)v_{(-1)}\ot v\0\ot h\1.
$$
There is a $C$-colinear map $V\ot H\to H\ot V$ such that 
$v\ot h\mapsto S(h)\ot v$. Since $S$ is bijective, this map is an isomorphism 
in $\CM$. Now Proposition 3.4 yields
$$
k1_C\sqC(H\ot V)\cong k1_C\sqC(V\ot H)\cong(H\ot k1_C)\sqC V.
$$
Since $H\ot k1_C\cong H$ in $\MC$, we finally get a $k$-linear bijection
$$
\Psi(H\ot V)\cong H\sqC V,\eqno(5.1)
$$
natural in $V$. If $\Psi$ is exact, then so is the functor $\Psi(H\ot{}?)$, 
whence $H\sqC{}?$ is exact on the category $\CM$ as well, establishing (a). 
Moreover, $H\sqC V=0$ if and only if $\Psi(H\ot V)=0$. If $\Psi$ is faithfully 
exact, the latter equality yields $V=0$, i.e., $H$ is right faithfully coflat 
over $C$.

$\mathrm{(b)}\Rar\mathrm{(c)}$. For each left $A$-module $W$ we have 
$\Phi(W)=H\ot_AW$ with the left $H$-module and $C$-comodule structures coming 
from those on $H$. The adjunction morphism
$$
\eta_W:W\to\Psi\Phi(W)=\lco C(H\ot_AW)\sbs H\ot_AW
$$
is the $A$-linear map sending $v\in W$ to $1\ot v\in H\ot_AW$. We have 
$\Phi(A)\cong H$. Up to this isomorphism the adjunction $\eta_A$ is the 
identity map $A\to \lco CH=A$. Since both $\Phi$ and $\Psi$ preserve direct 
sums, it follows that the adjunction $\eta_F$ is bijective whenever $F$ is a 
free $A$-module. Presenting $W$ as the cokernel of a homomorphism of free 
$A$-modules $F'\to F$, we get a commutative diagram
$$
\diagram{
F'&\hidewidth\lmapr5{}\hidewidth&F&\hidewidth\lmapr5{}\hidewidth&W
&\hidewidth\lmapr4{}&0\cr
\noalign{\smallskip}
\lmapdr{16}{\eta_{F'}}&&\lmapdr{16}{\eta_F}&&\lmapdr{16}{\eta_W}\cr
\noalign{\smallskip}
\Psi\Phi(F')&\lmapr2{}&\Psi\Phi(F)&\lmapr2{}&\Psi\Phi(W)&\lmapr2{}&0\cr
}
$$
in which the top row is exact and both $\eta_{F'}$ and $\eta_F$ are 
isomorphisms. The functor $\Phi$ is right exact. If $\Psi$ is exact, then the 
bottom row in the diagram is exact as well, and we deduce that $\eta_W$ is an 
isomorphism for each $W$. This means that $\Phi$ is fully faithful.

$\mathrm{(c)}\Rar\mathrm{(d)}$. If $\eta_W$ is an isomorphism, then the left 
$A$-module $W$ embeds in the induced $H$-module $H\ot_AW$.
\endproof

The proofs of Corollaries 5.1--5.3 dualize as well. We give below only 
statements of the respective conclusions.

\proclaim
Corollary 5.5.
Suppose that $K$ is a Hopf algebra, $H$ its factor Hopf algebra with 
bijective antipode, and $C$ a left $H$-module factor coalgebra of $H$. If 
$K$ is right (faithfully) coflat over $C,$ then so too is $H$.
\endproclaim

\proclaim
Corollary 5.6.
Suppose that $H$ is right coflat over its left module factor coalgebra $C$. 
Let $A=\lco CH$ and let $C'$ be any left $H$-module factor coalgebra of $H$ 
such that 
$$
A^+\sbs\Ker\pi'\sbs\Ker\pi
$$
where $\pi:H\to C$ and $\pi':H\to C'$ are the canonical surjections. Then 
$\lco CC'$ is spanned by the distinguished grouplike $\pi'(1)$ and $H$ is 
right coflat over $C'$.
\endproclaim

\proclaim
Corollary 5.7.
Suppose that $H$ is right coflat over its factor bialgebra $C$. Let 
$\pi$ be the canonical surjection $H\to C,$ and let $A=\lco CH$. Then:

\item(i)
$HA^+$ is the largest Hopf ideal of $H$ contained in $\Ker\pi$.

\item(ii)
$HA^+=S(HA^+)$ and $H$ is both left and right faithfully coflat over $H/HA^+$.

\item(iii)
$H$ is right faithfully coflat over $C$ if and only if $C$ is a Hopf 
factor algebra.

\endproclaim

\section
6. The generator property of Hopf modules

Bijectivity of the antipode $S:H\to H$ is not needed for the results of this 
section. Let $A$ be a right $H$-comodule algebra. We say that a Hopf module 
$M\in\MAH$ or $M\in\AMH$ is \emph{$A$-finite} if $M$ is finitely generated as 
an $A$-module. If $M$ and $N$ are two objects of $\AMH$ such that $M$ is 
$A$-finite, then $\Hom_A(M,N)$ is an $H$-comodule. This result of Caenepeel 
and Gu\'ed\'enon \cite{Cae-G04, Prop. 4.2} requires bijectivity of $S$, unlike 
its version for the category $\MAH$ presented below:

\proclaim
Proposition 6.1.
If $M$ is an $A$-finite object of $\MAH,$ then for any $N\in\MAH$ there is a 
right $H$-comodule structure on $\Hom_A(M,N)$ with respect to which the 
evaluation map $\,\Hom_A(M,N)\ot M\to N\,$ is $H$-colinear.
\endproclaim

\Proof.
Let $U$ be any finite-dimensional $H$-subcomodule of $M$ such that $M=UA$. 
The free $A$-module $F=U\ot A$ with the tensor product of $H$-comodule 
structures on $U$ and $A$ is an object of $\MAH$, and we get an exact sequence 
$0\to K\to F\to M\to0$ in $\MAH$. It gives rise to an exact sequence of vector 
spaces
$$
0\to\Hom_A(M,N)\to\Hom_A(F,N)\to\Hom_A(K,N).
$$
The dual $A$-module $\Hom_A(F,A)$ is generated freely by the vector space $U^*$ 
dual to $U$, whence a $k$-linear bijection
$$
\Hom_A(F,N)\cong N\ot_A\Hom_A(F,A)\cong N\ot_A(A\ot U^*)\cong N\ot U^*.
$$
There is an $H$-comodule structure which makes $U^*$ the left dual of $U$ in 
the monoidal category $\calM^H$. With the $H$-comodule structure on $\Hom_A(F,N)$ 
corresponding to the tensor product of $H$-comodule structures on $N$ and 
$U^*$ the evaluation map which factors as
$$
\Hom_A(F,N)\ot F\cong N\ot U^*\ot U\ot A\to N\ot k\ot A\cong N\ot A\to N
$$
is $H$-colinear since so are the evaluation map $U^*\ot U\to k$ and the module 
structure map $N\ot A\to N$. For each finite-dimensional $H$-subcomodule 
$V\sbs F$ the composite
$$
\diagram{
\Hom_A(F,N)\ot k\rlap{${}\cong\Hom_A(F,N)$}\cr
\downarrow&\cr
\Hom_A(F,N)\ot V\ot V^*&{}\hrar{}&\Hom_A(F,N)\ot F\ot V^*
\to N\ot V^*\cong\Hom_k(V,N)\cr
}
$$
obtained by means of the coevaluation $k\to V\ot V^*$ is $H$-colinear too, so 
it follows that the subspace $V^\flat\sbs\Hom_A(F,N)$ consisting of those 
$A$-module homomorphisms $F\to N$ that vanish on $V$ is stable under the 
coaction of $H$. Now $\Hom_A(M,N)$ is identified with 
$\,K^\flat=\{\xi\in\Hom_A(F,N)\mid\xi(K)=0\}$ which is the 
intersection of those subspaces $V^\flat$ when $V$ runs over all 
finite-dimensional $H$-subcomodules of $K$. Thus $\Hom_A(M,N)$ is an 
$H$-subcomodule of $\Hom_A(F,N)$.

The evaluation map $\,\Hom_A(M,N)\ot M\to N\,$ is $H$-colinear since so is its 
composite with the canonical surjection $\,\Hom_A(M,N)\ot F\to\Hom_A(M,N)\ot M$.
\endproof

Right $H$-comodules are identified with rational left modules for the dual 
algebra $H^*$ \cite{Sw}. In particular, the comodule $\Hom_A(M,N)$ in 
Proposition 6.1 is a rational left $H^*$-module.

In an arbitrary Hopf module $M\in\MAH$ the $A$-submodule $UA$ generated by any 
finite-dimensional $H$-subcomodule $U\sbs M$ is an $A$-finite subobject of $M$. 
It follows that each element of $M$ lies in an $A$-finite subobject. Hence 
the set $\calF(M)$ of $A$-finite subobjects of $M$ is directed by inclusion and
$$
\Hom_A(M,N)\cong\liminv_{M'\in\calF(M)}\Hom_A(M',N)
$$
in the category of vector spaces. It follows that $\Hom_A(M,N)$ has a 
uniquely determined left $H^*$-module structure such that the canonical maps 
$$
\Hom_A(M,N)\to\Hom_A(M',N)\eqno(6.1)
$$
are $H^*$-linear for all $M'\in\calF(M)$. This yields

\proclaim
Corollary 6.2.
For arbitrary $M,N\in\MAH$ there is a left $H^*$-module structure on 
$\Hom_A(M,N)$ which depends functorially on $M$ and $N$.
\endproclaim

\Remark.
The largest rational submodule $\HOM_A(M,N)$ of $\Hom_A(M,N)$ is the largest 
subspace of $\Hom_A(M,N)$ which has a right $H$-comodule structure such that 
the evaluation map $\,\HOM_A(M,N)\ot M\to N\,$ is $H$-colinear. This provides 
a different approach to the comodule structures introduced by Ulbrich 
\cite{Ulb90} which shows also that $\,\HOM_A(M,N)=\Hom_A(M,N)\,$ whenever $M$ 
is $A$-finite.
\endremark

The \emph{trace ideal} $\,T_M$ of a right $A$-module $M$ is a two-sided ideal 
of $A$ which coincides with the image of the evaluation map
$$
e_M:\ \Hom_A(M,A)\ot M\to A.\eqno(6.2)
$$
For $M$ to be a generator in $\MA$ it is necessary and sufficient that $T_M=A$.

\proclaim
Corollary 6.3.
If $M\in\MAH,$ then $T_M$ is stable under the coaction of $H$.
\endproclaim

\Proof.
Since $M$ is the union of its $A$-finite subobjects, it suffices to check that 
for each $A$-finite subobject $M'\sbs M$ the image of $\,\Hom_A(M,A)\ot M'$ 
under the evaluation map (6.2) is stable under the coaction of $H$. But the 
restriction of $e_M$ to $\,\Hom_A(M,A)\ot M'$ factors through the evaluation map
$$
e_{M'}:\ \Hom_A(M',A)\ot M'\to A,
$$
and therefore
$$
e_M\bigl(\Hom_A(M,A)\ot M'\bigr)=e_{M'}(I\ot M')
$$
where $I$ is the image of the canonical map (6.1) with $N=A$. Since (6.1) is 
an $H^*$-linear map, its image $I$ is an $H^*$-submodule, and therefore a 
subcomodule of the right $H$-comodule $\Hom_A(M',A)$. The desired conclusion 
follows immediately from the fact that $e_{M'}$ is a homomorphism of right 
$H$-comodules.
\endproof

\proclaim
Corollary 6.4.
Suppose that $A$ has no nonzero proper $H$-costable two-sided ideals. Then 
$M\in\MAH$ is a generator in $\MA$ if and only if $\,\Hom_A(M,A)\ne0$.
\endproclaim

\Proof.
By Corollary 6.3 and the current hypothesis the trace ideal $T_M$ has to be 
either 0 or the whole $A$. Thus $T_M=A$ if and only if $T_M\ne0$, but this is 
equivalent to the condition $\,\Hom_A(M,A)\ne0$.
\endproof

\proclaim
Corollary 6.5.
Suppose that $A$ is a right coideal subalgebra of $H$ such that $H$ is a 
generator in $\AM,$ and let $M\in\MAH$. Then $M$ is a generator in $\MA$ if 
and only if $\,\Hom_A(M,A)\ne0$.
\endproclaim

\Proof.
Suppose that $I$ is an $H$-costable two-sided ideal of $A$. Then $IH$ is a 
right ideal and a right coideal of $H$. Hence either $IH=0$ or $IH=H$. Since 
$H$ is a generator in $\AM$, we have $I=A\cap IH$, i.e., either $I=0$ or 
$I=A$. So we are in the situation of Corollary 6.4.
\endproof

\section
7. When does coflatness imply faithful coflatness?

Let $C$ be a left $H$-module factor coalgebra of $H$. We will be concerned 
with the case when $H$ is left coflat over $C$, i.e., the functor 
$?{}\sqCH:\MC\to\calM^H$ is exact. In this case $H$ is left faithfully coflat 
over $C$ if $V\sqCH\ne0$ for each nonzero right $C$-comodule $V$. By the local 
finiteness of comodules it suffices to verify this condition for simple 
comodules.

We will denote by $\Com_C$ and $\Com_H$ the vector spaces of $C$- and 
$H$-colinear maps, i.e., morphisms in the categories of comodules.

\setitemsize(b)
\proclaim
Lemma 7.1.
Let $V\in\MC$ and $U\in\calM^H$. Then $V\sqCH\ne0$ under any of the following 
two assumptions:

\item(a)
$\Com_C(U,V)\ne0,$

\item(b)
$\Com_C(V,U)\ne0$ and $H$ is left coflat over $C$.

\endproclaim

\Proof.
Since there is a canonical isomorphism of vector spaces 
$$
\Com_C(U,V)\cong\Com_H(U,\,V\sqCH),
$$
the conclusion in case (a) is immediately clear. In particular, 
$U\sqCH\ne0$ whenever $U$ is a nonzero right $H$-comodule.

Suppose now that $H$ is left coflat over $C$ and $\ph:V\to U$ is a nonzero 
$C$-colinear map. Then $V'=\Img\ph$ is a nonzero $C$-subcomodule of $U$. By 
exactness of the functor $?{}\sqCH$ the $H$-comodule $V'\sqCH$ is a 
homomorphic image of $V\sqCH$. So it suffices to show that $V'\sqCH\ne0$.

By the local finiteness of comodules we may assume that $\dim V'<\infty$ and 
also $\dim U<\infty$. In this case the dual vector space $U^*$ has a right 
$H$-comodule structure which makes $U^*$ the left dual of $U$ in the monoidal 
category $\calM^H$. With the trivial $H$-comodule structure on the base field $k$ 
the evaluation map $U^*\ot U\to k$ is a morphism in $\calM^H$, whence its 
restriction $\psi:U^*\ot V'\to k$ is a morphism in $\MC$. Here $U^*\ot V'$ is 
a $C$-comodule with respect to coaction (3.1). Since $\psi$ is obviously 
nonzero, it is surjective. But $k\sqCH\ne0$ since $k\in\calM^H\!$, 
and by exactness of the functor $?{}\sqCH$ it follows that
$$
U^*\ot(V'\sqCH)\cong(U^*\ot V')\sqCH\ne0
$$
where we have used isomorphism (3.4).
\endproof

CoFrobenius coalgebras were introduced by Lin \cite{Lin77} and the more 
general class of quasi-coFrobenius coalgebras by G\'omez Torrecillas and 
N\u{a}st\u{a}sescu \cite{Gom-N95}. For a thorough discussion of these 
notions we refer the reader to \cite{Das-NR}.

\proclaim
Lemma 7.2.
Suppose that $\dim C<\infty$. Then $C$ is left and right coFrobenius.
\endproclaim

\Proof.
The dual algebra $C^*$ embeds as a subalgebra in $H^*$. We have $C=H/I$ where 
$I$ is a left ideal of finite codimension in $H$. The annihilator $J$ of the 
left $H$-module $H/I$ is a two-sided ideal of finite codimension in $H$. Since 
$J\sbs I$, the image of $C^*$ in $H^*$ consists of linear functions vanishing 
on $J$, and so this image lies in the dual Hopf algebra $H^\circ$ of $H$. 
Since this image is stable under the right hit action of $H$ on $H^*$, it is 
in fact a finite-dimensional left coideal subalgebra of $H^\circ$. The algebra 
$H^*$ embeds in the direct product of the duals of finite-dimensional 
subcoalgebras of $H$. As a consequence $H^*$ is weakly finite. Hence so too is 
its subalgebra $H^\circ$. Applying \cite{Sk07, Th. 6.1}, we deduce that $C^*$ 
is a Frobenius algebra. This means that $C\cong C^*$ as $C^*$-modules on each 
side, and therefore $C$ is left and right coFrobenius.
\endproof

Note that Lemmas 7.1, 7.2 and items (a), (b), (c), (d) of the next proposition 
do not require bijectivity of the antipode.

\proclaim
Proposition 7.3.
Suppose that $H$ is left coflat over $C,$ and let $A=\lco CH$. Then $H$ is 
left faithfully coflat over $C$ in any of the following cases:

\item(a)
each simple $H$-comodule has dimension at most\/ $2,$

\item(b)
$\dim C<\infty,$

\item(c)
$C$ is either left or right quasi-coFrobenius,

\item(d)
$C$ is a factor Hopf algebra of $H,$

\item(e)
$A$ is right normal, i.e., stable under the right adjoint action of $H,$

\item(f)
$H$ is right flat over $A$.

\endproclaim

\Proof.
We have to prove that $V\sqCH\ne0$ for each simple right $C$-comodule $V$. 
This will follow from Lemma 7.1 once we show that such a comodule is isomorphic 
either to an $\MC$-subcomodule or to an $\MC$-factor comodule of some right 
$H$-comodule.

\smallskip
(a) Any simple right $C$-comodule $V$ occurs as a subfactor in some simple 
right $H$-comodule $U$. Since $\dim U\le2$, all $\MC$-composition series of 
$U$ have length 1 or 2, and the argument above does apply.

\smallskip
(b) This case is subsumed by case (c), in view of Lemma 7.2.

\smallskip
(c) If $C$ is left quasi-coFrobenius, then $C$ is projective in $\MC$ and a 
generator in $\CM$ by \cite{Gom-N95, Th. 1.3, Prop. 2.5}. More precisely, let 
$V\in\MC$ be a simple comodule. The injective hull $Q$ of $V$ is an 
indecomposable $\MC$-direct summand of $C$. Since $C$ is projective in $\MC$, 
so too is $Q$. Hence $\dim Q<\infty$ by \cite{Gom-N95, Lemma 1.2}. Then, by 
duality, $Q^*$ is a projective cover of $V^*$ in $\CM$, but also $Q^*$ is 
injective in $\CM$. It follows that $Q^*$ is a $\CM$-direct summand of $C$. As 
a consequence, there exist epimorphisms $H\to C\to Q^*\to V^*$ in $\CM$. 

Denoting by $f:H\to V^*$ the composite of these surjective $C$-colinear maps, 
we get $f(U)=V^*$ for some finite-dimensional left coideal $U$ of $H$. Then 
$V\cong V^{**}$ embeds as a $C$-subcomodule in the right $H$-comodule $U^*$.

By symmetry, if $C$ is right quasi-coFrobenius, then there exists an 
epimorphism $H\to V$ in $\MC$.

\smallskip
(d) In this case the conclusion was observed by Doi \cite{Doi83, Remark on p. 
247}.

\smallskip
(e) and (f)
\ Let $C=H/I$ where $I$ is a coideal and a left ideal of $H$. Put $I'=HA^+$ 
and $C'=H/I'$. By Lemma 1.1 $I'\sbs I$. Hence $C$ is a factor coalgebra of $C'$. 
Also, $\lco{C'}\!H=A$ (see Lemma 4.4). By Propositions 2.1 and 2.2 $H$ is a 
generator in $\AM$ and injective in $\CpM$. In other words, $H$ is left coflat 
over $C'$.

In case (e) $I'$ is a Hopf ideal by \cite{Tak94, Prop. 1.4}, and so $C'$ is a 
factor Hopf algebra of $H$. Therefore $H$ is left faithfully coflat over $C'$ 
by (d). In case (f) each finite-dimensional right $C'$-comodule embeds in some 
right $H$-comodule by Proposition 4.2. Then $H$ is left faithfully coflat over 
$C'$ by Lemma 7.1.

Thus it remains to prove that $C=C'$, i.e., $I=I'$. Note that $I$ and $I'$ are 
objects of $\HMC$. For any $M\in\HMC$ we have $M\sqCH\in\HMH$, and by the 
fundamental theorem on Hopf modules
$$
M\sqCH\cong H\ot(M\sqCH)\co H\cong H\ot M\co C.\eqno(7.1)
$$
Note that $C\co C=k1_C$. The canonical map $\pi:H\to C$ is a morphism in $\HMC$ 
with $I=\Ker\pi$. Put $B=H\co C$. By Lemma 1.1 $\pi(x)=\ep(x)1_C$ for $x\in B$. 
Hence 
$$
I\co C=I\cap B=B^+.
$$
By Lemma 1.6 $A=S(B)$, whence $I'=HB^+$ by Lemma 1.4. It follows then that 
$I'\,\vphantom{I}\co C=I\co C$, and by (7.1) the functor $?{}\sqCH$ takes the 
inclusion map $I'\to I$ to an isomorphism $I'\sqCH\cong I\sqCH$. Since that 
functor is exact, we deduce that
$$
(I/I')\sqCH=0.\eqno(7.2)
$$
Note that $I/I'$ is a right $C$-subcomodule of the coalgebra $C'=H/I'$. Since 
the exact functor $?{}\sqCH:\MC\to\calM^H$ is isomorphic to the composite of 
the functors
$$ 
?{}\sqC C':\MC\to\MCp\qquad\text{and}\qquad ?{}\sqCpH:\MCp\to\calM^H
$$
where the second one is faithfully exact, we conclude that the functor 
$?{}\sqC C'$ is exact as well, i.e., $C'$ is left coflat over $C$, and also 
$(I/I')\sqC C'=0$. In case (e) Lemma 7.1(b) applied to the Hopf algebra $C'$ 
and its left module factor coalgebra $C$ shows at once that $I/I'=0$.

In case (f) consider any simple right $C$-subcomodule $V\sbs C'$. There exists 
a finite-dimensional right coideal $W\sbs C'$ such that $V\sbs W$. By 
Proposition 4.2 $W$ embeds as a $C'$-subcomodule in some right $H$-comodule $U$. 
So it follows that $V$ is isomorphic to a $C$-subcomodule of $U$. But then 
$V\sqCH\ne0$ by Lemma 7.1, which implies $V\not\sbs I/I'$, in view of (7.2). 
In other words, $I/I'$ has zero socle as a right $C$-comodule. Again, 
this entails $I/I'=0$, and we are done.
\endproof

If $H$ is left coflat over $C$ and left flat over $A=\lco CH$, then $H$ is 
left faithfully coflat over $C'=H/HA^+$, but it is not clear whether $C=C'$, 
in contrast to case (f) of Proposition 7.3. This situation will reappear later 
in Corollary 10.7.

\proclaim
Theorem 7.4.
Suppose that $H$ is left coflat over $C$ and right flat over $A=\lco CH$. Then 
$H$ is left faithfully coflat over $C,$ left faithfully flat over $A$.

If, in addition, either $C$ has a proper right coideal of finite codimension 
or all simple left $A$-modules are finite dimensional, then $H$ is also right
faithfully coflat over $C,$ right faithfully flat over $A$.

\endproclaim

\Proof.
We have shown already in Proposition 7.3 that $H$ is left faithfully coflat over 
$C$. It is then a consequence of Theorem 1.7 that $C=H/HA^+$ and $H$ is left 
faithfully flat over $A$. It remains to prove the second part of the theorem.

By Takeuchi's Theorem 1 the functor $?{}\sqCH:\MC\to\MAH$ is an equivalence of 
categories. Thus for each $M\in\MAH$ there is $V\in\MC$ such that 
$M\cong V\sqCH$. If $\dim V<\infty$, then, by Proposition 4.2, there exists a 
right $H$-comodule $U$ in which $V$ embeds as a $C$-subcomodule. Then $M$ embeds 
in $U\sqCH\cong U\ot A$ where we have used Proposition 3.1. So $M$ embeds in a 
free $A$-module, which obviously implies that $\,\Hom_A(M,A)\ne0\,$ provided 
that $V\ne0$.

Suppose that $C$ has a right coideal $R\ne C$ of finite codimension. Put 
$V=C/R$. The canonical map $C\to V$ is an epimorphism in $\MC$. So it gives 
rise to an epimorphism $H\cong C\sqCH\to M=V\sqCH$ in $\MAH$. Since 
$\dim V<\infty$, the conclusion made in the preceding paragraph ensures the 
existence of a nonzero right $A$-linear map $H\to A$. But $H$ is a generator in 
$\AM$ by the Masuoka-Wigner theorem. Then $H$ is also a generator in $\MA$ by 
Corollary 6.5. Hence $H$ is right faithfully flat over $A$, and, by Theorem 
1.7, right faithfully coflat over $C$.

To deal with the other alternative in the hypothesis we will show that 
$H\ot_AV\ne0$ whenever $V$ is a nonzero finite-dimensional left $A$-module. 
This will be done by dualizing the argument in Lemma 7.1. The dual vector 
space $V^*$ is a right $A$-module. Consider the induced right $H$-module 
$T=V^*\ot_AH$. Since $H$ is left faithfully flat over $A$, the canonical 
$A$-linear map $V^*\to T$ is injective. It follows that the $A$-module 
$V\cong V^{**}$ is a homomorphic image of the left $H$-module $U=T^*$. In 
particular, there exists a nonzero $A$-linear map $U\to V$.

Since $\dim V<\infty$, there is a canonical $k$-linear bijection 
$\Hom_k(U,V)\cong V\ot U^*$. Here we equip $U^*$ with the left $H$-module 
structure making use of the antipode $S$. Then $V\ot U^*$ has the left 
$A$-module structure described in section 3. The corresponding $A$-module 
structure on $\Hom_k(U,V)$ is given by the formula
$$
(af)(u)=\sum a\1f\bigl(S(a\2)u\bigr),\qquad a\in A,\ f\in\Hom_k(U,V),\ u\in U.
$$
If some nonzero $f$ is $A$-linear, then $af=\ep(a)f$ for all $a\in A$, and so 
the 1-dimensional subspace spanned by $f$ is an $A$-submodule of $\Hom_k(U,V)$ 
isomorphic to $A/A^+$. This shows that the right $A$-module $A/A^+$ embeds in 
$V\ot U^*$. Applying the exact functor $H\ot_A{}?$, we embed $H/HA^+$ in 
$H\ot_A(V\ot U^*)$. Since $H\ne HA^+$, it follows that
$$
H\ot_A(V\ot U^*)\ne0.
$$
But $H\ot_A(V\ot U^*)\cong(H\ot_AV)\ot U^*$ by Proposition 3.5. Hence 
$H\ot_AV\ne0$, as claimed.

If all simple left $A$-modules are finite dimensional, we conclude that 
$H\ot_AV\ne0$ for each simple left $A$-module $V$. This implies again that the 
flat right $A$-module $H$ is faithfully flat.
\endproof

\proclaim
Corollary 7.5.
Suppose that $C$ is left quasi-coFrobenius. If $H$ is left coflat over 
$C,$ then $H$ is both left and right faithfully coflat over $C$.
\endproclaim

\Proof.
Put $A=\lco CH$. By \cite{Gom-N95, Prop. 2.5} $C$ is projective in $\MC$. The 
canonical surjection $H\to C$ is then a split epimorphism in $\MC$. Hence $H$ 
has an $\MC$-direct summand isomorphic to $C$. It follows that $H$ is a 
cogenerator in $\MC$, but then $H$ is projective, and therefore flat, in $\MA$ 
by Proposition 2.2. Also, each indecomposable $\MC$-direct summand $Q$ of $C$ 
is projective in $\MC$, whence $\dim Q<\infty$ by \cite{Gom-N95, Lemma 1.2}. 
This implies that $C=Q\oplus I$ where $I$ is a right coideal of finite 
codimension in $C$. Thus Theorem 7.4 applies.
\endproof

\proclaim
Corollary 7.6.
Suppose that $A$ is a finite-dimensional quasi-Frobenius right coideal 
subalgebra of $H$. Then $H$ is left flat over $A$ if and only if $H$ is right 
flat over $A$. In this case $H$ is faithfully flat over $A$ on both sides.
\endproclaim

\Proof.
Take $C=H/HA^+$. Since $A$ is selfinjective on both sides, it splits off in 
$H$ as an $\MA$-direct summand, and also as an $\AM$-direct summand. By 
Proposition 2.1 $H$ is left and right coflat over $C$, while $A=\lco CH$ by 
Corollary 4.5. The desired conclusion follows from Theorem 7.4.
\endproof

In the next corollary bijectivity of the antipode of $H$ is deduced from the 
other assumptions about $H$ (see \cite{Sk21}).

\proclaim
Corollary 7.7.
Suppose that a Hopf algebra $H$ as an ordinary algebra is finitely generated, 
noetherian and satisfies a polynomial identity. If $H$ is either left or right 
coflat over its left module factor coalgebra $C,$ then $H$ is both left and 
right faithfully coflat over $C$.
\endproclaim

\Proof.
By \cite{Sk21, Cor. 4.6} $H$ is left and right flat over any right coideal 
subalgebra. In particular, this holds for $A=\lco CH$. Theorem 7.4 applies 
again (in the case of right coflatness we pass to the Hopf algebra $H\cop$).

Since $H$ satisfies a polynomial identity, so too does $A$. By Kaplansky's 
theorem each primitive factor ring of $A$ is simple artinian (see, e.g., 
\cite{Pro, Ch. 2, Th. 1.1}). In particular, the annihilators of simple 
$A$-modules are maximal ideals of $A$. Since $H$ is left faithfully flat over 
$A$, each right $A$-module embeds in the induced $H$-module. It follows that 
each simple right $A$-module embeds in a simple right $H$-module. All simple 
$H$-modules are finite dimensional over $k$ since $H$ is a finitely generated 
PI algebra \cite{Pro, Ch. 5, Th. 1.2}. Hence so too are all simple right 
$A$-modules. It follows that all maximal ideals of $A$ have finite 
codimension, and therefore all simple left $A$-modules are also finite 
dimensional. This validates the use of the second part of Theorem 7.4.
\endproof

\section
8. Projectivity and injectivity of Hopf modules

This section deals with properties (0.1)--(0.4) of right coideal subalgebras 
and left module factor coalgebras. We continue to assume that $H$ is a Hopf 
algebra with bijective antipode.  Let $A$ be a right coideal subalgebra and 
$C$ a left $H$-module factor coalgebra of $H$. Denote by $\pi$ the canonical 
surjection $H\to C$.

\setitemsize(ii)
\proclaim
Lemma 8.1.
Suppose that all objects of $\MAH$ are flat $A$-modules. Then:

\item(i)
All nonzero objects of $\MAH$ are projective generators in $\MA$.

\item(ii)
$H$ is left and right faithfully flat over $A$.

\endproclaim

\Proof.
Since $H/A$ is an object of $\MAH$, it is right flat over $A$, and it follows 
that $H$ is right faithfully flat over $A$.

It is well-known that finitely presented flat modules are projective. Thus any 
Hopf module $M\in\MAH$ which is finitely presented as an $A$-module is 
projective in $\MA$. By Corollary 6.3 the trace ideal $T_M$ of such a module 
is stable under the coaction of $H$. Then $HT_M$ is a left ideal and a right 
coideal of $H$. If $M\ne0$, then $T_M\ne0$. It follows then that $HT_M=H$ since 
$H$ is a simple object of the category $\HMH$, whence $T_M=A$ by the right 
faithful flatness of $H$ over $A$. This means that $M$ is a projective 
generator in $\MA$, and so, in particular, $M\ne MA^+$.

If $M\in\MAH$ is finitely presented in $\MA$ and $M'$ is any $A$-finite 
subobject of $M$, then $M/M'$ is also finitely presented in $\MA$, whence 
$M/M'$ is projective in $\MA$ as we have just seen. In this case the exact 
sequence $0\to M'\to M\to M/M'\to0$ splits in $\MA$. So $M'$ is an $A$-module 
direct summand of $M$, which implies that $M'$ is finitely presented in $\MA$ 
as well and $M'/M'A^+$ embeds in $M/MA^+$. Moreover, if $M'\ne M$, then $M/M'$ 
is a projective generator in $\MA$, whence $M'/M'A^+\ne M/MA^+$.

In other words, the $(H,A)$-Hopf modules finitely presented in $\MA$ form an 
abelian full subcategory of $\MAH$, and the functor $M\mapsto M/MA^+$ is 
faithfully exact on this subcategory. In particular, when $M$ is finitely 
presented in $\MA$, the lattice of $A$-finite subobjects of $M$ embeds into 
the lattice of subspaces of the finite-dimensional vector space $M/MA^+$. 
We deduce then that the $A$-finite subobjects of $M$ satisfy the ascending 
chain condition. Since every object of $\MAH$ is the sum of its $A$-finite 
subobjects, it follows that all subobjects of $M$ have to be $A$-finite. 
This shows that each factor object of $M$ is finitely presented 
in $\MA$ whenever so is $M$.

Suppose now that $M$ is an arbitrary object of $\MAH$. If $M'\sbs M$ is any 
$A$-finite subobject, then $M'=UA$ for some finite-dimensional $H$-subcomodule  
$U\sbs M$. There is also a Hopf module structure on $U\ot A$, and the action 
of $A$ on $M$ yields an epimorphism $U\ot A\to M'$ in $\MAH$. Since $U\ot A$ 
is a free $A$-module of finite rank, we deduce that $M'$ is finitely presented, 
and therefore projective, in $\MA$. The set of $A$-finite subobjects of $M$ is 
directed by inclusion, and the union of this set is the whole $M$. This 
enables one to prove by applying Zorn's lemma that for an epimorphism $V\to W$ 
in $\MA$ any $A$-linear map $M\to W$ can be lifted to an $A$-linear map $M\to 
V$. Hence $M$ is projective in $\MA$. Each subobject $M'\sbs M$ is an 
$\MA$-direct summand since $M/M'$ is also projective in $\MA$. We know already 
that nonzero $A$-finite subobjects are generators in $\MA$. Hence so too is 
$M$ provided that $M\ne0$.

Thus (i) has been verified. It follows that every short exact sequence in 
$\MAH$ splits in $\MA$. But then Takeuchi's functor $\Phi:\MAH\to\MC$ is 
faithfully exact, and finally Lemma 4.1 ensures that $H$ is left faithfully 
flat over $A$.
\endproof

\proclaim
Lemma 8.2.
Suppose that all objects of $\HCM$ are coflat $C$-comodules. Then:

\item(i)
All nonzero objects of $\HCM$ are injective cogenerators in $\CM$.

\item(ii)
$H$ is left and right faithfully coflat over $C$.

\endproclaim

\Proof.
Coflat comodules are injective. Hence all objects of $\HCM$ are injective in 
$\CM$, and therefore exact sequences in $\HCM$ split in $\CM$. In particular, 
the canonical surjection $\pi:H\to C$ is a split epimorphism in $\CM$. This 
implies that $H$ is left faithfully coflat over $C$. We also see that 
Takeuchi's functor $\HCM\to\AM$ is exact.

Let $M\in\HCM$. Then $M$ is coflat over $C$. We will prove faithful coflatness 
provided that $M\ne0$, which is equivalent to the property of being an 
injective cogenerator in $\CM$. It suffices to show that $V\sqCM\ne0$ for 
each nonzero right coideal $V$ of $C$. Since the functor $?{}\sqCM$ commutes 
with directed colimits there exists a largest right coideal $R$ of $C$ 
such that $R\,\sqCM=0$. (In fact $R$ is a subcoalgebra since $R$ is stable 
under all $\MC$-endomorphisms of $C$.) 

We have to show that $R=0$. By formula (3.7) $(U\ot R)\sqCM=0$ for each 
right $H$-comodule $U$. In particular, this holds for $U=H$. Since the module 
structure map $H\ot C\to C$ is a morphism in $\MC$ and the functor $?{}\sqCM$ 
is exact, we get $(HR)\sqCM=0$. Hence $HR=R$ by the maximality of $R$, i.e., 
$R$ is an $H$-submodule of $C$. Put
$$
I=\{h\in H\mid(\pi\ot\id)\De(h)\in R\ot H\}\cong R\sqCH.
$$
Then $I$ is a left ideal and a right coideal of $H$. Since $H$ is a simple 
object of the category $\HMH$, it follows that either $I=0$ or $I=H$. Note 
that
$$
H/I\cong(C/R)\sqCH.
$$
By the left faithful coflatness of $H$ over $C$ we must have either $R=0$ or 
$R=C$. But $R\ne C$ since $C\sqCM\cong M\ne0$. Hence $R=0$, as claimed.

Faithful coflatness of $M$ implies that $\lco CM\ne0$, which means that 
Takeuchi's functor $\HCM\to\AM$ is also faithful. Hence $H$ is right 
faithfully coflat over $C$ by Lemma 5.4.
\endproof

\proclaim
Lemma 8.3.
Suppose that $A\sbs\lco CH$ and either $A$ satisfies $(0.1)$ or $C$ 
satisfies $(0.2)$. Then there is a natural isomorphism of vector spaces
$$
(V\sqCH)\ot_AW\cong V\sqC(H\ot_AW)\eqno(8.1)
$$
for each right $C$-comodule $V$ and left $A$-module $W$.
\endproclaim

\Proof.
The left $C$-comodule structure on $H$ is given by the map
$$
\la=(\pi\ot\id)\circ\De:H\to C\ot H.
$$
Since $A\sbs\lco CH$, the action of $A$ on $H$ by right multiplications commutes 
with the coaction $\la$. The cotensor product $V\sqCH$ is the vector space 
defined by means of the exact sequence
$$
0\mapr{}V\sqCH\mapr{}V\ot H\lmapr8{\rho\ot\id-\id\ot\la}V\ot C\ot H\eqno(8.2)
$$
where $\rho:V\to V\ot C$ is the right $C$-comodule structure on $V$. We view 
$H$ as an object of $\MAH$ with the action of $A$ by right multiplications and 
the right coaction of $H$ by its comultiplication. Then $V\ot H$ and 
$V\ot C\ot H$ are objects of $\MAH$ with respect to the action of $A$ and the 
coaction of $H$ defined on $H$. This makes (8.2) an exact sequence in $\MAH$ 
and $V\sqCH$ a subobject of $V\ot H$.

Likewise there is an exact sequence in $\HCM$
$$
H\ot A\ot W\lmapr8{\mu\ot\id-\id\ot\nu}H\ot W\mapr{}H\ot_AW\mapr{}0\eqno(8.3)
$$
where $\mu:H\ot A\to H$ and $\nu:A\ot W\to W$ are, respectively, the restriction 
of the multiplication in $H$ and the module structure map. Here $H$ acts by left 
multiplications on itself and $C$ coacts on $H$ by means of $\la$.

If $A$ satisfies (0.1), then all monomorphisms and epimorphisms of the 
category $\MAH$ split in $\MA$. It follows that exactness of (8.2) is preserved 
under the functor $?{}\ot_AW$, whence $(V\sqCH)\ot_AW$ is isomorphic to the 
kernel of the map which defines the right hand side of (8.1). Similarly, if 
$C$ satisfies (0.2), then exactness of (8.3) is preserved under the functor 
$V\sqC{}?$, which again yields (8.1).
\endproof

\proclaim
Theorem 8.4.
Takeuchi's correspondence gives a bijection between the set of right coideal 
subalgebras of $H$ satisfying condition $(0.1)$ and the set of left $H$-module 
factor coalgebras of $H$ satisfying condition $(0.2)$.

There is a similar bijection with respect to conditions $(0.3)$ and $(0.4)$.
\endproclaim

\Proof.
Suppose that $A$ is a right coideal subalgebra of $H$ satisfying condition 
(0.1), and let $C=H/HA^+$. By Lemma 8.1 $H$ is left and right faithfully flat 
over $A$. Then $A=\lco CH$ and $H$ is both left and right faithfully coflat 
over $C$ by Takeuchi's Theorem 1 and its generalization in Theorem 1.7.

Let $M\in\HCM$. By Takeuchi's Theorem 2 there exists a left $A$-module $W$ such 
that $M\cong H\ot_AW$. By Lemma 8.3 the functor $?{}\sqC M:\CM\to\calM_k$ is 
isomorphic to the composite of the functor $?{}\sqC H:\CM\to\MAH$, which is 
exact by coflatness of $H$ over $C$, and the functor 
$?{}\ot_AW:\MAH\to\calM_k$, which is exact because all exact sequences in 
$\MAH$ split in $\MA$. Thus $M$ is coflat, and therefore injective, in $\CM$. 
This shows that $C$ satisfies (0.2).

Conversely, suppose that $C$ is any left $H$-module factor coalgebra of $H$ 
satisfying (0.2), and let $A=\lco CH$. By Lemma 8.2 $H$ is left and right 
faithfully coflat over $C$. Then $C=H/HA^+$ and $H$ is left and right 
faithfully flat over $A$.
Each Hopf module $M\in\MAH$ is isomorphic to $V\sqCH$ for some right 
$C$-comodule $V$. The corresponding functor $M\ot_A{}?:\AM\to\calM_k$ is 
isomorphic to the composite of the functors
$$
H\ot_A{}?:\AM\to\HCM\qquad\text{and}\qquad V\sqC{}?:\HCM\to\calM_k,
$$
both of which are exact. Thus each module $M\in\MAH$ is flat in $\MA$, and 
Lemma 8.1 shows that $A$ satisfies (0.1).

Consider now the second pair of conditions on $A$ and $C$. We may identify 
$\AMH$ with $\calM_{A\op}^{H\op}$ and $\HMC$ with $\HCopM$. Suppose that $A$ 
is a right coideal subalgebra of $H$ satisfying (0.3). Then the already proved 
part of Theorem 8.4 applied to the Hopf algebra $H\op$ and its right coideal 
subalgebra $A\op$ satisfying (0.1) ensures that the factor coalgebra 
$D=H/A^+H$ has the property that all objects of the category $\DMH=\HopDM$ are 
injective in $\DM$. It follows from the category equivalences of Lemma 1.8 
that $C=H/HA^+$ satisfies (0.4).

Conversely, suppose that $C$ is a left $H$-module factor coalgebra of $H$ 
satisfying condition (0.4). Then $C\cop$ is a left $H\cop$-module 
factor coalgebra of the Hopf algebra $H\cop$ satisfying (0.2), whence 
all objects of the category $\HMB=\MB^{H\cop}$ where $B=H\co C$ are projective 
in $\MB$. Applying Lemma 1.9, we deduce that $A=\lco CH$ satisfies (0.3).
\endproof

\proclaim
Corollary 8.5.
Suppose that $H$ is a Hopf algebra with cocommutative coradical. Then for any 
left $H$-module factor coalgebra $C$ of $H$ all nonzero Hopf modules in the 
categories $\HCM$ and $\HMC$ are injective cogenerators in $\CM$ and 
$\MC\!,$ respectively.
\endproclaim

\Proof.
Let $A=\lco CH$. By \cite{Ma91, Cor. 1.5} all Hopf modules in the categories 
$\MAH$ and $\AMH$ are projective $A$-modules. Hence Theorem 8.4 and Lemma 8.2 
apply.
\endproof

\proclaim
Theorem 8.6.
Suppose that $H$ is left and right faithfully flat over its right coideal 
subalgebra $A,$ and let $C=H/HA^+$. If either $A$ is an algebra of finite weak 
global dimension or $C$ is a coalgebra of finite global dimension, then
properties $(0.1)$--$(0.4)$ are all satisfied.
\endproclaim

\Proof.
By Takeuchi's Theorem 1 and Theorem 1.7 $H$ is left and right faithfully coflat 
over $C$, and $A=\lco CH$. Moreover, the functor $?{}\sqCH$ gives an equivalence 
of categories $\MC\to\MAH$. This functor takes cofree comodules to direct sums 
of copies of $H$ which are projective right $A$-modules by the Masuoka-Wigner 
theorem. Any injective right $C$-comodule $Q$ is a direct summand of a cofree 
comodule, whence $Q\sqCH$ is projective in $\MA$. If $V\in\MC$ has a finite 
injective resolution

$$
0\to V\to Q_0\to Q_1\to\cdots\to Q_n\to0,
$$
then, applying the functor $?{}\sqCH$, we get an exact sequence
$$
0\to V\sqCH\to Q_0\sqCH\to Q_1\sqCH\to\cdots\to Q_n\sqCH\to0
$$
in $\MAH$. This sequence splits in $\MA$ since each $Q_i\sqCH$ is projective 
in $\MA$. Hence $V\sqCH$ is an $\MA$-direct summand of $Q_0\sqCH$, and 
therefore $V\sqCH$ is projective in $\MA$.

Suppose that $C$ is a coalgebra of finite global dimension. Since each object 
of $\MAH$ is isomorphic to $V\sqCH$ for some comodule $V\in\MC$ which has a 
finite injective dimension, we deduce that (0.1) is satisfied. Applying this 
fact to the Hopf algebra $H\cop$ with $A$ replaced by $B=H\co C=S^{-1}(A)$ and 
$C$ replaced by $C\cop$, we also conclude that all objects of the category 
$\HMB$ are projective in $\MB$. Lemma 1.9 shows that (0.3) is satisfied as 
well. Then (0.2) and (0.4) hold by Theorem 8.4.

In the other thread of the proof we exploit the functor $H\ot_A{}?$ which gives 
an equivalence of categories $\AM\to\HCM$ by Takeuchi's Theorem 2. This functor 
takes free $A$-modules to direct sums of copies of $H$ which are coflat in 
$\CM$. Any flat left $A$-module $F$ is a directed colimit of free modules by 
Lazard's theorem. Since tensor products commute with directed colimits, 
$H\ot_AF$ is a directed colimit of coflat $C$-comodules, and therefore 
$H\ot_AF$ is itself coflat, i.e., injective in $\CM$. If a left $A$-module $V$ 
has a finite flat resolution
$$
0\to F_n\to\cdots\to F_1\to F_0\to V\to0,
$$
then there is an exact sequence
$$
0\to H\ot_AF_n\to\cdots\to H\ot_AF_1\to H\ot_AF_0\to H\ot_AV\to0,
$$
in $\HCM$ which splits in $\CM$ since each $H\ot_AF_i$ is injective in $\CM$. 
In this case $H\ot_AV$ is a $\CM$-direct summand of $H\ot_AF_0$, and 
therefore injective in $\CM$.

Suppose that $A$ is an algebra of finite weak global dimension. Then the 
preceding argument verifies (0.2). Switching to the Hopf algebra $H\op$ with 
$A$ replaced by $A\op$ and $C$ replaced by $D=H/A^+H$, we also conclude that 
all objects of the category $\DMH$ are injective in $\DM$. By Lemma 1.8 
this yields (0.4). Finally, (0.1) and (0.3) follow from Theorem 8.4.
\endproof

\section
9. Vanishing of Ext and Tor

Although we are not able to prove the conclusion of Theorem 8.6 when finiteness 
of global dimensions is not assumed, we will show in Proposition 9.4 that Hopf 
modules still satisfy somewhat weaker properties. This result is preceded by 
three lemmas used in its proof.

Let us view $H$ as a left $H$-module and a right $H$-comodule in the usual way. 
For each left $A$-module $V$ the tensor product $V\ot H$ over the base field 
is an object of the category $\AMH$ with respect to the action of $A$ given by 
formula (3.8) and the coaction of $H$ defined on $H$ (see \cite{Doi92, 1.4}).

\proclaim
Lemma 9.1.
The functor $\AM\to\AMH$ defined by the assignment $V\mapsto V\ot H$ is right 
adjoint of the forgetful functor $\AMH\to\AM$. In other words, there are 
natural bijections
$$
\AMH(M,\,V\ot H)\cong\Hom_A(M,V)
$$
for objects $M\in\AMH$ and $V\in\AM$. As a consequence, $Q\ot H$ is injective 
in $\AMH$ whenever $Q$ is an injective left $A$-module.
\endproclaim

\Proof.
The adjointness of these functors is a special case of \cite{Br99, Cor. 3.7}. 
There is a canonical bijection between the $H$-colinear maps $M\to V\ot H$ 
and the $k$-linear maps $M\to V$ given by the assignment 
$\ph\mapsto(\id\ot\ep)\circ\ph$. Under this bijection $\ph$ is $A$-linear if 
and only if so is the corresponding map $M\to V$. Since the forgetful functor 
$\AMH\to\AM$ is exact, its right adjoint preserves injectives.
\endproof

\proclaim
Lemma 9.2.
Let $\Phi:\AMH\to\MD$ and $\Phi':\AM\to\HCM$ be Takeuchi's functors where 
$C=H/HA^+$ and $D=H/A^+H$. Denote by $i_S$ the composite of the forgetful 
functor $\HCM\to\CM$ and the equivalence $\CM\to\MD$ arising from the 
antiisomorphism of coalgebras $C\to D$ induced by the antipode $S$. 
There are natural isomorphisms
$$
\Phi(V\ot H)\cong(i_S\circ\Phi')(V)\quad\text{for $\,V\in\AM$}.
$$
\endproclaim

\Proof.
Recall that $\Phi(M)=M/A^+M$ for $M\in\AMH$. Let $M=V\ot H$. For $v\in V$ and 
$h\in H$ put $\overline{v\ot h}=v\ot h+A^+M\in\Phi(M)$. Similarly to the case 
of the category $\MAH$ considered in \cite{Sk10, Lemma 1.2} there is a 
$k$-linear bijective map
$$
\xi:\Phi(M)\to H\ot_AV=\Phi'(V)
$$
such that $\,\xi(\overline{v\ot h})=S^{-1}h\ot v\,$ and 
$\xi^{-1}(h\ot v)=\overline{v\ot Sh}$ for $v\in V$ and $h\in H$.
Let $\pi:H\to C$ and $\pi':H\to D$ be the canonical 
surjections. The right $D$-comodule structure on $\Phi(M)$ is 
$$
\overline{v\ot h}\mapsto\sum\,(\overline{v\ot h\1})\ot\pi'(h\2).
$$
It corresponds under $\xi$ to the right $D$-comodule structure on $H\ot_AV$ 
given by the assignment
$$
h\ot v\mapsto\sum\,(h\2\ot v)\ot\pi'(Sh\1),
$$
while the left $C$-comodule structure on $\Phi'(V)$ is given by
$$
h\ot v\mapsto\sum\,\pi(h\1)\ot(h\2\ot v).
$$
The antiisomorphism of coalgebras $C\to D$ induced by $S$ sends $\pi(h)$ to
$\pi'(Sh)$ for each $h\in H$.
\endproof

Dually, the forgetful functor $\HMC\to\MC$ for a left $H$-module factor 
coalgebra $C$ of $H$ has a left adjoint $\MC\to\HMC$ defined by the assignment 
$V\mapsto H\ot V$ with the $C$-comodule structure given by formula (3.1). 
Similarly to Lemma 9.2 we have

\proclaim
Lemma 9.3.
Let $\Psi:\HMC\to\BM$ and $\Psi':\MC\to\MAH$ be Takeuchi's functors where 
$A=\lco CH$ and $B=H\co C$. Denote by $i_S$ the composite of the forgetful 
functor $\MAH\to\MA$ and the equivalence $\MA\to\BM$ arising from the 
antiisomorphism of algebras $B\to A$ induced by the antipode $S$. 
There are natural isomorphisms
$$
\Psi(H\ot V)\cong(i_S\circ\Psi')(V)\quad\text{for $\,V\in\MC$}.
$$
\endproclaim

\Proof.
Bijection (5.1) with $H$ replaced by $H\cop$ gives an isomorphism of vector 
spaces
$$
\Psi(H\ot V)=(H\ot V)\co C\cong V\sqCH=\Psi'(V).
$$
Here $(H\ot V)\co C$ is a $B$-submodule of $H\ot V$ and $V\sqCH$ is an 
$A$-submodule of $V\ot H$ with respect to the left and right actions defined 
by the rules
$$
b\cdot(h\ot v)=bh\ot v,\qquad(v\ot h)\cdot a=v\ot ha
$$
for $a\in A$, $b\in B$, $h\in H$, $v\in V$. The above isomorphism is induced 
by the $k$-linear map $H\ot V\to V\ot H$ such that $h\ot v\mapsto v\ot Sh$. 
This map is in fact $B$-linear, assuming that $V\ot H$ is made into a left 
$B$-module by means of the functor $i_S$. Hence the requested isomorphism in 
$\BM$.
\endproof

For the lack of a better alternative we will denote by $\Ext_C^i$ the $\Ext$ 
groups for the abelian category $\MC$. The notation $\Cotor_C^i$ will be used 
for the right derived functors of the cotensor product $\sqC$. The latter 
functors were introduced in a paper by Eilenberg and Moore \cite{Ei-M65}.

\setitemsize(iii)
\proclaim
Proposition 9.4.
Suppose that $H$ is left and right faithfully flat over $A$. Let $V$ be a left 
$A$-module. If either $V$ has finite flat dimension or the $A$-module structure 
on $V$ extends to an $H$-module structure on $V,$ then:

\item(i)
$\,V\ot H$ is injective in $\AMH$ and projective in $\AM,$

\item(ii)
$\,\Ext_A^i(M,\,V)=0$ for all Hopf modules $M\in\AMH$ and $i>0,$

\item(iii)
$\,\Tor^A_i(M,\,V)=0$ for all Hopf modules $M\in\MAH$ and $i>0$.

Put $C=H/HA^+$. Given a right $C$-comodule $W$ such that either $W$ has 
finite injective dimension or the $C$-comodule structure $W\to W\ot C$ lifts 
to an $H$-comodule structure $W\to W\ot H,$ we have:

\item(iv)
$\,H\ot W$ is projective in $\HMC$ and injective in $\MC,$

\item(v)
$\,\Ext_C^i(W,\,M)=0$ for all Hopf modules $M\in\HMC$ and $i>0,$

\item(vi)
$\,\Cotor_C^i(W,\,M)=0$ for all Hopf modules $M\in\HCM$ and $i>0$.

\endproclaim

\Proof.
(i) Let $\Phi'$ and $\Phi$ be as in Lemma 9.2. Suppose that $V$ is a left 
$A$-module of finite flat dimension. In this case it was shown in the proof of 
Theorem 8.6 that $\Phi'(V)=H\ot_AV$ is injective in $\CM$. Then $\Phi(V\ot H)$ 
is injective in $\MD$ by Lemma 9.2. Hence $V\ot H$ is injective in $\AMH$ 
since $\Phi$ is an equivalence by Corollary 1.2. Let $V\cong F/K$ where $F$ 
is a free $A$-module. Then the flat dimension of the $A$-module $K$ is finite 
as well, and therefore $K\ot H$ is injective in $\AMH$ as we have just shown. 
It follows that the exact sequence
$$
0\to K\ot H\to F\ot H\to V\ot H\to0
$$
splits in $\AMH$. In other words, $V\ot H$ is a direct summand of $F\ot H$. 
But $F\ot H$ is a free $A$-module by Corollary 3.7. Hence $V\ot H$ is 
projective in $\AM$.

Suppose now that $V$ is a left $H$-module. Then $V\ot H$ is an object of 
$\HMH$. By the fundamental theorem on Hopf modules in this category $V\ot H$ 
is isomorphic to a direct sum of copies of $H$. Since $H$ is projective in 
$\AM$ by the Masuoka-Wigner theorem, so too is $V\ot H$. Since $\Phi(H)=D$ is 
injective in $\MD$, so is the $D$-comodule $\Phi(V\ot H)$, whence $V\ot H$ is 
injective in $\AMH$.

(ii) For any $M\in\AMH$ there exists an epimorphism $F_0\to M$ in $\AMH$ with 
$F_0$ being $A$-free. Indeed, it suffices to take $F_0=A\ot U$ where $U$ is 
any $H$-subcomodule of $M$ such that $M=AU$. The map $F_0\to M$ is given by 
the action of $A$ on $M$. Iterating, we get a resolution $\cdots\to F_1\to 
F_0\to M\to0$ in $\AMH$ such that each $F_i$ is a free $A$-module. Since 
$V\ot H$ is injective in $\AMH$ by (i), the induced sequence of vector spaces
$$
0\to\AMH(M,\,V\ot H)\to\AMH(F_0,\,V\ot H)\to\AMH(F_1,\,V\ot H)\to\cdots
$$
is exact. By Lemma 9.1 this sequence is isomorphic to the sequence
$$
0\to\Hom_A(M,V)\to\Hom_A(F_0,V)\to\Hom_A(F_1,V)\to\cdots
$$
whose cohomology computes the groups $\,\Ext_A^i(M,\,V)$.

(iii) Let $M\in\MAH$. As in (ii) we can find a resolution $\cdots\to F_1\to 
F_0\to M\to0$ in $\MAH$ such that each $F_i$ is a free $A$-module. It gives 
rise to an exact sequence of vector spaces
$$
\cdots\to F_1\ot_A(V\ot H)\to F_0\ot_A(V\ot H)\to M\ot_A(V\ot H)\to0
$$
since $V\ot H$ is projective in $\AM$. By Proposition 3.8 this means 
that the sequence
$$
\cdots\to(F_1\ot_AV)\ot H\to(F_0\ot_AV)\ot H\to(M\ot_AV)\ot H\to0
$$
is exact. Hence so is the sequence 
$\cdots\to F_1\ot_AV\to F_0\ot_AV\to M\ot_AV\to0$ whose homology computes 
the groups $\,\Tor^A_i(M,\,V)$.

\smallskip
The dual assertions (iv)--(vi) are proved similarly with the aid of Lemma 9.3.
\endproof

\section
10. Quotient categories and localization

Recall that the quotient category $\calA/\calS$ of an abelian category $\calA$ 
is defined when $\calS$ is a \emph{Serre subcategory}, i.e., a full subcategory 
such that for each exact sequence
$$
0\to X'\to X\to X''\to0\quad\text{in $\calA$}
$$
one has $X\in\calS$ if and only if $X'\in\calS$ and $X''\in\calS$ 
(see \cite{Gab62, Ch. III} and \cite{Fai, Ch. 15}). The category $\calA/\calS$ 
has the same objects as $\calA$, and there is an exact canonical functor 
$\calA\to\calA/\calS$ which is the identity map on objects and which sends an 
arbitrary morphism $f$ in $\calA$ to an isomorphism in $\calA/\calS$ if and 
only if both $\Ker f$ and $\Coker f$ are in $\calS$.

A Serre subcategory $\calS$ is \emph{localizing} (respectively, 
\emph{colocalizing}) if the canonical functor $\calA\to\calA/\calS$ admits a 
right (respectively, left) adjoint $\calA/\calS\to\calA$ called the 
\emph{section functor}. The latter is always fully faithful. The composite 
$\calA\to\calA/\calS\to\calA$ of the canonical functor and the section functor 
is called the \emph{(co)localization functor}.
The localizing subcategories are interchanged with the colocalizing ones when 
passing to the dual abelian category.

Quotient categories arise generally as follows. If an exact functor 
$T:\calA\to\calB$ between two abelian categories admits a fully faithful right 
(respectively, left) adjoint $S:\calB\to\calA$, then $\Ker T$ is 
a (co)localizing subcategory of $\calA$ and $T$ induces an equivalence of 
categories $\calA/\Ker T\approx\calB$ (see \cite{Gab62, Ch. III, Prop. 5}). 
The composite $ST$ is then the respective (co)localization functor.

In the case of a Grothendieck category $\calA$ its Serre subcategory $\calS$ is 
localizing if and only if each object of $\calA$ has a largest subobject lying 
in $\calS$. For example, given a ring $A$ and a flat left $A$-module $M$ the 
right modules $V$ such that $V\ot_AM=0$ form a localizing subcategory of the 
category $\MA$ of right $A$-modules. Similarly, for a coflat left comodule $M$ 
over a coalgebra $C$ there is a localizing subcategory of $\MC$ consisting of 
right comodules $V$ such that $V\sqCM=0$.

Every colocalizing subcategory of a Grothendieck category is also a localizing 
subcategory \cite{Nas-T96, Prop. 2.2}. In $\MC$ the localizing subcategories
are precisely the categories $\calM^U$ of comodules over subcoalgebras $U\sbs C$ 
such that $U=U\wedge U$ with respect to the wedge product. 
N\u{a}st\u{a}sescu and Torrecillas found a necessary and sufficient condition 
for such a subcategory to be colocalizing \cite{Nas-T96, Prop. 3.1}.

\smallskip
We are going to explain how quotient categories arise in Hopf theory. 
Retaining our standing assumption about the Hopf algebra $H$, let $A$ be a 
right coideal subalgebra of $H$. Consider full subcategories $\TA\sbs\MA$ and 
$\TAH\sbs\MAH$ consisting of modules characterized by the property that each 
element is annihilated by some nonzero $H$-costable right ideal of $A$. Define 
similarly subcategories $\AT\sbs\AM$ and $\ATH\sbs\AMH$ in terms of 
$H$-costable left ideals.

\proclaim
Lemma 10.1.
Suppose that $H$ is right flat over $A$. Given a left $A$-module $V$, we have 
$\,V\in\AT\,$ if and only if $\,H\ot_AV=0$.
\endproclaim

\Proof.
If $I$ is any left ideal of $A$ stable under the coaction $H$, then $HI$ is a 
left ideal and a right coideal of $H$. It follows that $HI=H$ provided that 
$I\ne0$. If $v\in V$ is annihilated by such a left ideal $I$, then $1\ot v=0$ 
in the $H$-module $H\ot_AV$. This shows that $H\ot_AV=0$ whenever $V\in\AT$.

Conversely, suppose that $H\ot_AV=0$. Equip $M=V\ot H$ with the left action of 
$A$ defined by formula (3.8) and the right coaction of $H$ given by the map 
$\id\ot\De$. Then $M\in\AMH$. Since $H\in\MAH$, we have
$$
H\ot_AM\cong(H\ot_AV)\ot H=0\eqno(10.1)
$$
by Proposition 3.8. We may also view $A$ as an object of $\AMH$. The category 
$\AMH$ is a right module category over the monoidal category $\calM^H$ (cf. 
discussion in section 3). Each finite-dimensional $H$-subcomodule $U\sbs M$ 
has a left dual $U^*$ in $\calM^H$, whence a canonical bijection
$$
\AMH(A\ot U,\,M)\cong\AMH(A,\,M\ot U^*).
$$
The action of $A$ on $M$ defines a morphism $A\ot U\to M$ in $\AMH$. It 
corresponds to a morphism $f:A\to M\ot U^*$ such that $\Ker f$ coincides with 
the annihilator of $U$ in $A$. Thus $I=\Ker f$ is an $H$-costable left ideal 
of $A$ which annihilates all elements of $U$. Furthermore, $f$ induces an 
embedding of $A/I$ in $M\ot U^*$. Since $H\ot_A{}?$ is an exact functor which 
annihilates $M\ot U^*$ by (10.1), we get $H\ot_AA/I=0$, whence $I\ne0$. Since 
each element of $M$ is contained in a finite-dimensional $H$-subcomodule, it 
follows that $M\in\AT$. There is an $A$-linear surjection $M\to V$ given by 
the map $\id_V\ot\ep$. Hence $V\in\AT$ too.
\endproof

\setitemsize(iii)
\proclaim
Theorem 10.2.
Put $\,C=H/HA^+$. If $H$ is left flat over $A,$ then:

\item(i)
$\TA$ and $\TAH$ are localizing subcategories, respectively, of $\MA$ and $\MAH,$

\item(ii)
there are category equivalences $\,\MAH/\TAH\approx\MC$ 
and $\,\MA/\TA\approx\HMC$.

\smallskip
If $H$ is right flat over $A,$ then:

\item(iii)
$\AT$ and $\ATH$ are localizing subcategories, respectively, of $\AM$ and $\AMH,$

\item(iv)
there are category equivalences $\,\AMH/\ATH\approx\CM\,$ 
and $\,\AM/\AT\approx\HCM$.

\endproclaim

\Proof.
The cases of left and right flatness are interchanged with each other by 
moving from $H$ to $H\op$. Suppose that $H$ is left flat over $A$. By Lemma 
10.1 applied to the Hopf algebra $H\op$ the objects of $\TA$ are all right 
$A$-modules $V$ such that $V\ot_AH=0$. By the flatness assumption this defines 
a localizing subcategory of $\MA$. We see also that $\TAH$ is a localizing 
subcategory of $\MAH$ consisting of Hopf modules $M$ such that $M\ot_AH=0$.

Let $\Phi$ be Takeuchi's functor $\MAH\to\MC$ and $\Psi$ its right adjoint.
By Lemma 4.1 $\Phi$ is exact and $\Psi$ is fully faithful.
$\Psi$. This implies that $\Ker\Phi$ is a localizing subcategory of $\MAH$ and 
$\Phi$ induces an equivalence of the quotient category $\MAH/\Ker\Phi$ with 
the category $\MC$ (see \cite{Gab62, Ch. III, Prop. 5}). The kernel of $\Phi$ 
consists of all objects $M\in\MAH$ such that $M=MA^+$. It follows from (2.1) 
that $\,\Ker\Phi=\TAH$. This yields the first equivalence in (ii).

The second equivalence in (ii) will follow from (iv). If the second 
equivalence in (iv) holds for the Hopf algebra $H\op$, then 
$\,\MA/\TA\approx\HopDM\approx\DMH\approx\HMC\,$ where $D=H/A^+H$ (see Lemma 1.8).

So let us verify the second equivalence in (iv) assuming now that $H$ is right 
flat over $A$. Let $\Phi$ be Takeuchi's functor $\AM\to\HCM$ and $\Psi$ its 
right adjoint. Since $\Phi=H\ot_A{}?$, this functor is exact and $\Ker\Phi=\AT$. 
It remains to prove that $\Psi$ is fully faithful. This amounts to showing that 
the adjunction morphisms
$$
\eta_M:\Phi\Psi(M)=H\ot_A\lco CM\to M,\qquad M\in\HCM,
$$
are isomorphisms. Note that $\eta_M$ is given by the action of $H$ on $M$. It 
follows already from Propositions 4.2 and 4.3 that $\eta_M$ is surjective for 
each $M$.

The forgetful functor $\HCM\to\HM$ has a right adjoint $\HM\to\HCM$ given by 
the assignment $U\mapsto C\ot U$ (cf. Lemma 9.1). The Hopf module structure on 
$C\ot U$ is obtained from the structure available on $C$ by applying the 
tensoring functor $?{}\ot U$, as discussed in section 3. Since 
$\lco CC=k1_C\cong A/A^+$, the map $\eta_C$ is bijective. It follows then from 
Proposition 3.5 that $\eta_{\,C\ot U}$ is bijective for any $U\in\HM$.

Given an arbitrary object $M\in\HCM$, take $U$ to be any left $H$-module such 
that there exists a monomorphism $f:M\to U$ in $\HM$. Let $g:M\to C\ot U$ be 
the morphism in $\HCM$ corresponding to $f$. Then $(\ep\ot\id_U)\circ g=f$, 
whence $g$ is injective. Setting $N=C\ot U$, we get a commutative diagram
$$
\diagram{
H\ot_A\lco CM & \lmapr3{\eta_M} & M \cr
\noalign{\smallskip}
\lmapd{16}{\id\ot g'}{} && \lmapd{16}{}g \cr
\noalign{\smallskip}
H\ot_A\lco CN & \lmapr3{\eta_N} & N \cr
}
$$
where $g':\lco CM\to\lco CN$ is the restriction of the map $g$. Since 
$\,\id\ot g'$ and $\eta_N$ are both injective, so is $\eta_M$. Hence $\eta_M$ 
is bijective for each $M$, and we are done.
\endproof

\proclaim
Corollary 10.3.
Denote by $\calF$ the set of left ideals of $A$ such that a left ideal $I$ is in 
$\calF$ if and only if $I$ contains a nonzero $H$-costable left ideal of $A$. If 
$H$ is right flat over $A$ then $\calF$ is a left Gabriel topology.
\endproclaim

\Proof.
For each left ideal $I$ of $A$ we have $I\in\calF$ if and only if $A/I\in\AT$. 
In other words, $\calF$ is the filter of left ideals corresponding to the 
localizing subcategory $\AT$ (see \cite{Fai}, \cite{Gab62}, \cite{St}).
\endproof

\proclaim
Corollary 10.4.
If $H$ is right flat over $A$ and right coflat over $C=H/HA^+,$ then the algebra 
$A'=\lco CH$ is a perfect left localization of $A$. In other words, $A'$ is flat 
in $\MA$ and the inclusion map $A\to A'$ is an epimorphism in the category of 
rings.
\endproclaim

\Proof.
Under present assumptions $\Phi:\AM\to\HCM$ and $\Psi:\HCM\to\AM$ in (1.2) are 
both exact. Indeed, the functor $\Psi$ factors through the category $\ApM$, 
and so its exactness was already observed in Lemma 5.4. The localization 
functor $\Psi\Phi$ with respect to the localizing subcategory $\AT$ of $\AM$ 
is thus exact. It also preserves arbitrary direct sums. The algebra 
$A'\cong\Psi\Phi(A)$ is the localization of its subalgebra $A$. It is a 
general fact that in such a situation $A'$ is a perfect left localization of 
$A$ (see condition (c) in \cite{St, Ch. XI, Prop. 3.4}).

Right flatness of $A'$ over $A$ can be seen as follows. There are natural maps
$$
\xi_V:A'\ot_AV\to\lco C(H\ot_AV)=\Psi\Phi(V)
$$
which are bijective for all $V\in\AM$. Indeed, $\xi_A$ is bijective since 
$\Phi(A)\cong H$ and $\Psi(H)=A'$. Hence $\xi_V$ is bijective whenever $V$ is 
a free $A$-module. For arbitrary $V$ we take a presentation $F_1\to F_0\to 
V\to0$ where $F_0$ and $F_1$ are free $A$-modules and use right exactness of 
all functors. Since $\Psi\Phi$ is exact, so is the functor $\,A'\ot_A{}?$.

Finally, since $A'/A\in\AT$, this $A$-module is annihilated by $\Phi$, and 
therefore by the functor $\,A'\ot_A{}?$. Hence the canonical map 
$A'\ot_AA'\to A'$ is bijective, which is a necessary and sufficient condition 
for the ring extension $A\to A'$ to be epimorphic \cite{St, Ch. XI, Prop. 1.2}.
\endproof

Now we dualize previous results. Let $C$ be a left $H$-module factor coalgebra 
of the Hopf algebra $H$ and $\pi:H\to C$ the canonical surjection. 

Denote by $\TC$ the full subcategory of $\MC$ consisting of right comodules $V$ 
characterized by the property that for each element $v\in V$ there exists some 
$H$-stable subcoalgebra $R\sbs C$ such that $R\ne C$ and $\rho(v)\in V\ot R$ 
where $\rho:V\to V\ot C$ is the comodule structure map. Let $\HTC$ be the full 
subcategory of $\HMC$ consisting of Hopf modules which are in $\CT$ with 
forgotten $H$-module structure. The categories $\CT$ and $\HCT$ are defined 
similarly. For each subcoalgebra $R\sbs C$ the category of comodules $\calM^R$ 
may be identified with a full subcategory of $\MC$.

\proclaim
Lemma 10.5.
Suppose that $H$ is left coflat over $C$. Then there exists a largest proper 
$H$-stable subcoalgebra $R\sbs C$. For a right $C$-comodule $V$ we\/ have 
$V\in\TC$ if and only if\/ $V\in\calM^R,$ if and only if\/ $V\sqCH=0$.
\endproclaim

\Proof.
If $R$ is any $H$-stable subcoalgebra of $C$, then the set
$$
I=\{h\in H\mid(\pi\ot\id)\De(h)\in R\ot H\}\cong R\sqCH
$$
is a left ideal and a right coideal of $H$, and we have 
$\pi(h)=(\pi\ot\ep)\De(h)\in R$ for all $h\in I$, i.e., $\pi(I)\sbs R$. If 
$R\ne C$, then $I\ne H$, whence $I=0$ since $H$ is a simple 
object of the category $\HMH$, and therefore $\,R\sqCH=0$.

Now let $R$ be the sum of all proper $H$-stable subcoalgebras of $C$. So $R$ is 
itself an $H$-stable subcoalgebra and contains every $H$-stable subcoalgebra of 
$C$ other than the whole $C$. Since $?{}\sqCH$ is an exact functor commuting 
with directed colimits, we get $R\sqCH=0$. Hence $R\ne C$, yielding the 
first assertion of the lemma. It is clear thereof that $V\in\TC$ if and only 
if $\rho(V)\sbs V\ot R$, i.e., $V\in\calM^R$.

Furthermore, the coalgebra $C$ has a largest right coideal $J$ such that 
$J\,\sqCH=0$. For each $\MC$-endomorphism $\ph$ of $C$ we have $\ph(J)\sbs J$ 
since $\ph(J)\sqCH=0$ by exactness of $?{}\sqCH$. Hence $J$ is a subcoalgebra 
of $C$. By Proposition 3.4
$$
(H\ot J)\sqCH=0
$$
where the coaction of $C$ on $H\ot J$ is defined by formula (3.1). Since the 
action of $H$ gives a morphism $H\ot C\to C$ in $\MC$, we get $(HJ)\sqCH=0$. 
Hence $HJ=J$ by maximality of $J$, i.e. $J$ is stable under the action of 
$H$. It follows that $R=J$.
 
There exists a $C$-colinear embedding of $V$ in the direct sum of some family 
of copies of $C$. By exactness of the functor $?{}\sqCH$ the equality 
$f(V)\sqCH=0$ holds for a $C$-colinear map $f:V\to C$ whenever $V\sqCH=0$. 
Considering projections onto direct summands we deduce that $V\sqCH=0$ if and 
only if $f(V)\sqCH=0$ for all $f\in\Com_C(V,C)$. This is further equivalent to 
the condition that $f(V)\sbs R$ for all $f$, and therefore to the condition 
that $V$ embeds in a direct sum of copies of $R$, i.e, $V\in\calM^R$.
\endproof

\proclaim
Theorem 10.6.
Put $A=\lco CH$. If $H$ is left coflat over $C,$ then:

\item(i)
$\TC$ and $\HTC$ are colocalizing subcategories, respectively, of $\MC$ and 
$\HMC,$

\item(ii)
there are category equivalences $\,\HMC/\HTC\approx\MA$ 
and $\,\MC/\TC\approx\MAH$.

\smallskip
If $H$ is right coflat over $C,$ then:

\item(iii)
$\CT$ and $\HCT$ are colocalizing subcategories, respectively, of 
$\,\CM$ and $\HCM,$

\item(iv)
there are category equivalences $\,\HCM/\HCT\approx\AM\,$ 
and $\,\CM/\,\CT\approx\AMH$.

\endproclaim

\Proof.
First suppose that $H$ is right coflat over $C$. By Lemma 5.4 Takeuchi's 
functor $\,\Psi:\HCM\to\AM\,$ is exact and its left adjoint 
$\,\Phi:\AM\to\HCM\,$ is fully faithful. Hence $\Ker\Psi$ is a colocalizing 
subcategory of $\HCM$ and $\Psi$ induces an equivalence 
$\,\HCM/\Ker\Psi\approx\AM$. Recall that $\Psi(M)=\lco CM$ for $M\in\HCM$. We 
see from (2.2) that $\Psi(M)=0$ if and only if $H\sqCM=0$. Lemma 10.5 applied 
to the Hopf algebra $H\cop$ shows that $\,\Ker\Psi=\HCT$.

Suppose that $H$ is left coflat over $C$. Then Takeuchi's functor 
$\,\Psi:\MC\to\MAH$ is exact. Moreover, $\Ker\Psi=\TC$ by Lemma 10.5. Recalling 
the left adjoint functor $\Phi$ in (1.1), let us show that the adjunctions
$$
\xi_M:M\to(M/MA^+)\sqCH=\Psi\Phi(M),\qquad M\in\MAH,
$$
are isomorphisms, and so $\Phi$ is fully faithful. First, $\xi_A$ is bijective 
since $\Phi(A)\cong k1_C$ and $\Psi(k1_C)\cong\lco CH=A$. By Proposition 3.1 
$\xi_{U\ot A}$ is then bijective for any right $H$-comodule $U$. For 
an arbitrary $M\in\MAH$ we can find in $\MAH$ an exact sequence
$$
U_1\ot A\to U_0\ot A\to M\to0
$$ 
for some $U_0,\,U_1\in\calM^H$. Since both $\Phi$ and $\Psi$ are right exact, 
the induced sequence
$$
\Psi\Phi(U_1\ot A)\to\Psi\Phi(U_0\ot A)\to\Psi\Phi(M)\to0
$$
is exact, and it follows from naturality of $\xi$ that $\xi_M$ is bijective, 
yielding our claim. It implies again that $\Ker\Psi$ is a colocalizing 
subcategory of $\MC$ and $\Psi$ induces an equivalence 
$\,\MC/\Ker\Psi\approx\MAH$.
The remaining assertions of Theorem 10.6 follow from the already proved 
ones with $H$ replaced by $H\cop$.
\endproof

\proclaim
Corollary 10.7.
If $H$ is left coflat over $C$ and left flat over $A=\lco CH,$ then the 
coalgebra $C'=H/HA^+$ is a perfect right colocalization of $C$. So 
$C'$ is coflat in $\CM$ and the canonical surjection $C'\to C$ is a
monomorphism in the category of coalgebras.
\endproclaim

\Proof.
Both functors in (1.1) are exact and $C'\cong\Phi\Psi(C)$ is the 
colocalization of $C$ with respect to the colocalizing subcategory $\TC$ of 
$\MC$. Thus we are in the situation of \cite{Nas-T96, Th. 4.3}.
\endproof

It can only be repeated that there are no known examples of factor coalgebras 
$C$ for which the hypothesis of Corollary 10.7 is satisfied with $C'\ne C$ 
since the equality $C'=C$ does hold whenever $H$ is left faithfully coflat 
over $C$.

In contrast to this remark nontrivial examples of perfect localizations 
provided by earlier Corollary 10.4 are abundant already in the case of 
commutative Hopf algebras. Let $H=k[G]$ be the algebra of regular functions on 
an affine algebraic group $G$. Suppose that $G$ acts on an affine algebraic 
variety $X$ with a dense open orbit $O\sbs X$. Let $G_x$ be the stabilizer of 
a point $x\in O$. Then the quotient $G/G_x$ is a quasiaffine variety isomorphic 
to $O$. The algebra $A=k[X]$ embeds $G$-equivariantly in $A'=k[O]$, and the 
latter is isomorphic to the subalgebra of $G_x$-invariant elements of $H$ with 
respect to the action of $G_x$ by left translations. Thus $A$ and $A'$ are 
identified with right coideal subalgebras of $H$. Furthermore, $H$ is coflat 
over its Hopf factor algebra $C=H/HA^+\cong k[G_x]$ if and only if the orbit $O$ 
is an affine variety. In this case we meet the hypothesis of Corollary 10.4, 
and $A'\ne A$ provided that $O\ne X$.

Needless to say, Corollaries 10.3, 10.4, and 10.7 are valid also with the left 
and right sides interchanged.

\references
\nextref Bia-HM63
\auth{A.,Bia{\l}ynicki-Birula;G.,Hochschild;G.D.,Mostow}
\paper{Extensions of representations of algebraic linear groups}
\journal{Amer. J. Math.}
\Vol{85}
\Year{1963}
\Pages{131-144}

\nextref Br99
\auth{T.,Brzezi\'nski}
\paper{On modules associated to coalgebra Galois extensions}
\journal{J.~Algebra}
\Vol{215}
\Year{1999}
\Pages{290-317}

\nextref Br-W
\auth{T.,Brzezi\'nski;R.,Wisbauer}
\book{Corings and Comodules}
\publisher{Cambridge Univ. Press}
\Year{2003}

\nextref Cae-G04
\auth{S.,Caenepeel;T.,Gu\'ed\'enon}
\paper{Projectivity of a relative Hopf module over the subring of coinvariants}
\InBook{Hopf algebras}
\publisher{Marcel Dekker}
\Year{2004}
\Pages{97-108}

\nextref Cl-PS77
\auth{E.,Cline;B.,Parshall;L.,Scott}
\paper{Induced modules and affine quotients}
\journal{Math. Ann.}
\Vol{230}
\Year{1977}
\Pages{1-14}

\nextref Das-NR
\auth{S.,D\u asc\u alescu;C.,N\u ast\u asescu;S.,Raianu}
\book{Hopf Algebras, an Introduction}
\publisher{Marcel Dekker}
\Year{2000}

\nextref Doi83
\auth{Y.,Doi}
\paper{On the structure of relative Hopf modules}
\journal{Comm. Algebra}
\Vol{11}
\Year{1983}
\Pages{243-255}

\nextref Doi92
\auth{Y.,Doi}
\paper{Unifying Hopf modules}
\journal{J.~Algebra}
\Vol{153}
\Year{1992}
\Pages{373-385}

\nextref Ei-M65
\auth{S.,Eilenberg;J.C.,Moore}
\paper{Homology and fibrations I: Coalgebras, cotensor product and its derived functors}
\journal{Comment. Math. Helv.}
\Vol{40}
\Year{1965}
\Pages{199-236}

\nextref Fai
\auth{C.,Faith}
\book{Algebra: Rings, Modules and Categories I}
\publisher{Springer}
\Year{1973}

\nextref Gab62
\auth{P.,Gabriel}
\paper{Des cat\'egories ab\'eliennes}
\journal{Bull. Soc. Math. France}
\Vol{90}
\Year{1962}
\Pages{323-448}

\nextref Gom-N95
\auth{J.,G\'omez Torrecillas;C.,N\u ast\u asescu}
\paper{Quasi-co-Frobenius coalgebras}
\journal{J.~Algebra}
\Vol{174}
\Year{1995}
\Pages{909-923}

\nextref Gross
\auth{F.D.,Grosshans}
\book{Algebraic Homogeneous Spaces and Invariant Theory}
\BkSer{Lecture Notes Math.}
\BkVol{1673}
\publisher{Springer}
\Year{1997}

\nextref Ja
\auth{J.,Jantzen}
\book{Representations of Algebraic Groups}
\publisher{Academic Press}
\Year{1987}

\nextref Kop93
\auth{M.,Koppinen}
\paper{Coideal subalgebras in Hopf algebras: Freeness, integrals, smash products}
\journal{Comm. Algebra}
\Vol{21}
\Year{1993}
\Pages{427-444}

\nextref Lin77
\auth{B.I-p.,Lin}
\paper{Semiperfect coalgebras}
\journal{J.~Algebra}
\Vol{49}
\Year{1977}
\Pages{357-373}

\nextref Ma91
\auth{A.,Masuoka}
\paper{On Hopf algebras with cocommutative coradicals}
\journal{J.~Algebra}
\Vol{144}
\Year{1991}
\Pages{451-466}

\nextref Ma94
\auth{A.,Masuoka}
\paper{Quotient theory of Hopf algebras}
\InBook{Advances in Hopf algebras}
\publisher{Marcel Dekker}
\Year{1994}
\Pages{107-133}

\nextref Ma-D92
\auth{A.,Masuoka;Y.,Doi}
\paper{Generalization of cleft comodule algebras}
\journal{Comm. Algebra}
\Vol{20}
\Year{1992}
\Pages{3703-3721}

\nextref Ma-W94
\auth{A.,Masuoka;D.,Wigner}
\paper{Faithful flatness of Hopf algebras}
\journal{J.~Algebra}
\Vol{170}
\Year{1994}
\Pages{156-164}

\nextref Nas-T96
\auth{C.,N\u ast\u asescu;B.,Torrecillas}
\paper{Colocalization on Grothendieck categories with applications to coalgebras}
\journal{J.~Algebra}
\Vol{185}
\Year{1996}
\Pages{108-124}

\nextref Nas-TZ96
\auth{C.,N\u ast\u asescu;B.,Torrecillas;Y.H.,Zhang}
\paper{Hereditary coalgebras}
\journal{Comm. Algebra}
\Vol{24}
\Year{1996}
\Pages{1521-1528}

\nextref Ob77
\auth{U.,Oberst}
\paper{Affine Quotientenschemata nach affinen, algebraischen Gruppen und induzierte Darstellungen}
\journal{J.~Algebra}
\Vol{44}
\Year{1977}
\Pages{503-538}

\nextref Par-W91
\auth{B.,Parshall;J.-p.,Wang}
\paper{Quantum linear groups}
\journal{Mem. Amer. Math. Soc.}
\Vol{439}
\Year{1991}
\Pages{1-157}

\nextref Pro
\auth{C.,Procesi}
\book{Rings with Polynomial Identities}
\publisher{Marcel Dekker}
\Year{1973}

\nextref Scha00
\auth{P.,Schauenburg}
\paper{Faithful flatness over Hopf subalgebras: counterexamples}
\InBook{Interactions between Ring Theory and Representations of Algebras}
\publisher{Marcel Dekker}
\Year{2000}
\Pages{331-344}

\nextref Scha-Schn05
\auth{P.,Schauenburg;H.-J.,Schneider}
\paper{On generalized Hopf Galois extensions}
\journal{J.~Pure Appl. Algebra}
\Vol{202}
\Year{2005}
\Pages{168-194}

\nextref Schn92
\auth{H.-J.,Schneider}
\paper{Normal basis and transitivity of crossed products for Hopf algebras}
\journal{J.~Algebra}
\Vol{152}
\Year{1992}
\Pages{289-312}

\nextref Schn93
\auth{H.-J.,Schneider}
\paper{Some remarks on exact sequences of quantum groups}
\journal{Comm. Algebra}
\Vol{21}
\Year{1993}
\Pages{3337-3357}

\nextref Scho
\auth{A.H.,Schofield}
\book{Representation of Rings over Skew Fields}
\publisher{Cambridge Univ. Press}
\Year{1985}

\nextref Sk07
\auth{S.,Skryabin}
\paper{Projectivity and freeness over comodule algebras}
\journal{Trans. Amer. Math. Soc.}
\Vol{359}
\Year{2007}
\Pages{2597-2623}

\nextref Sk10
\auth{S.,Skryabin}
\paper{Models of quasiprojective homogeneous spaces for Hopf algebras}
\journal{J.~Reine Angew. Math.}
\Vol{643}
\Year{2010}
\Pages{201-236}

\nextref Sk21
\auth{S.,Skryabin}
\paper{Flatness over PI coideal subalgebras}
\journal{Israel J. Math.}
\Vol{245}
\Year{2021}
\Pages{735-772}

\nextref St
\auth{B.,Stenstr\"om}
\book{Rings of Quotients}
\publisher{Springer}
\Year{1975}

\nextref Sw
\auth{M.E.,Sweedler}
\book{Hopf Algebras}
\publisher{Benjamin}
\Year{1969}

\nextref Tak77
\auth{M.,Takeuchi}
\paper{Formal schemes over fields}
\journal{Comm. Algebra}
\Vol{5}
\Year{1977}
\Pages{1483-1528}

\nextref Tak79
\auth{M.,Takeuchi}
\paper{Relative Hopf modules---equivalences and freeness criteria}
\journal{J.~Algebra}
\Vol{60}
\Year{1979}
\Pages{452-471}

\nextref Tak94
\auth{M.,Takeuchi}
\paper{Quotient spaces for Hopf algebras}
\journal{Comm. Algebra}
\Vol{22}
\Year{1994}
\Pages{2503-2523}

\nextref Ulb90
\auth{K.-H.,Ulbrich}
\paper{Smash products and comodules of linear maps}
\journal{Tsukuba J. Math.}
\Vol{14}
\Year{1990}
\Pages{371-378}

\endreferences
\bye